\documentclass[a4paper,french]{article}

\usepackage[utf8]{inputenc}

\usepackage[french]{babel}

\usepackage{amsmath}\usepackage{epsf,amsfonts,amsthm}

\usepackage{authblk}  

\usepackage[hmargin=2.2cm,vmargin=2.5cm]{geometry}		
\newcommand{\be}{\begin{equation}}
\newcommand{\ee}{\end{equation}}
\newcommand{\bea}{\begin{eqnarray*}}
\newcommand{\eea}{\end{eqnarray*}}

\newcommand{\N}{\mathbb{N}}

\mathchardef\za="710B  
\mathchardef\zb="710C  
\mathchardef\zg="710D  
\mathchardef\zd="710E  
\mathchardef\zve="710F 
\mathchardef\zz="7110  
\mathchardef\zh="7111  
\mathchardef\zy="7112 
\mathchardef\zi="7113  
\mathchardef\zk="7114  
\mathchardef\zl="7115  
\mathchardef\zm="7116  
\mathchardef\zn="7117  
\mathchardef\zx="7118  
\mathchardef\zp="7119  
\mathchardef\zr="711A  
\mathchardef\zs="711B  
\mathchardef\zt="711C  
\mathchardef\zu="711D  
\mathchardef\zf="711E 
\mathchardef\zq="711F  
\mathchardef\zc="7120  
\mathchardef\zw="7121  
\mathchardef\ze="7122  
\mathchardef\zvy="7123  
\mathchardef\zvw="7124  
\mathchardef\zvr="7125 
\mathchardef\zvs="7126 
\mathchardef\zvf="7127  
\mathchardef\zG="7000  
\mathchardef\zD="7001  
\mathchardef\zY="7002  
\mathchardef\zL="7003  
\mathchardef\zX="7004  
\mathchardef\zP="7005  
\mathchardef\zS="7006  
\mathchardef\zU="7007  
\mathchardef\zF="7008  
\mathchardef\zW="700A  

\newtheorem{rem}{Remarque}
\newtheorem{theo}{Th\'eor\`eme}
\newtheorem*{theo*}{Th\'eor\`eme}
\newtheorem{prop}{Proposition}
\newtheorem{lem}{Lemme}

\newtheorem{defi}{D\'efinition}

\newenvironment{preuve}{\paragraph{D\'emonstration~:}}{\hfill$\square$}

\begin{document}

\title{Infinis morphismes de Leibniz pour les crochets d\'eriv\'es}

\author[1]{ Camille Laurent-Gengoux \thanks{camille.laurent-gengoux@univ-lorraine.fr}}
\author[1]{Mohsen Masmoudi  \thanks{mohsen.masmoudi@univ-lorraine.fr}}
\affil[1]{{\small Institut Elie Cartan de Lorraine (IECL), UMR 7502,  Universit\'e de Lorraine, Metz et Nancy, France.}}
 \maketitle

\begin{abstract}
\noindent
{\bf{Version fran\c{c}aise.}}
Il est connu que le crochet d\'eriv\'e d'un \'el\'ement de Maurer-Cartan d'une alg\`ebre de Lie diff\'erentielle gradu\'ee (DGLA) d\'efinit une alg\`ebre de Leibniz diff\'erentielle gradu\'ee.
Il est connu aussi que un morphisme de Lie infini entre DGLAs envoie un \'el\'ement de Maurer-Cartan sur un autre \'el\'ement de Maurer-Cartan. 
\'Etant donn\'es un morphisme Lie-infini,  un \'el\'ement de Maurer-Cartan et son image, nous construisons entre leurs alg\`ebres de Leibniz diff\'erentielles gradu\'ees
un morphisme de Leibniz infini, et ce de fa\c{c}on totalement explicite. Nous utilisons cette construction pour retrouver une formule de Dominique Manchon \`a propos du commutateur du produit-\'etoile.

\noindent
{\bf{English Abstract.}} The derived bracket of a Maurer-Cartan element in a differential graded Lie algebra (DGLA) is well-known to define a differential graded Leibniz algebra.
It is also well-known that a Lie infinity morphism between DGLAs maps a Maurer-Cartan element to a Maurer-Cartan element. 
Given a Lie-infinity morphism, a Maurer-element and its image, we show that both derived differential graded Leibniz algebras are related by a Leibniz-infinity morphism, 
and we construct it explicitely.
As an application, we recover a well-known formula of Dominique Manchon about the commutator of the star-product.

\noindent
{\bf Mots-clefs}: Alg\`ebres de Leibniz, alg\`ebre de Lie-infinies, formalit\'e et quantification. 

\end{abstract}

\tableofcontents

\section{Introduction}

Il y a dans la litt\'erature sur la quantification deux faits bien connus concernant les alg\`ebres de Lie gradu\'ees diff\'erentielles  (que l'on appellera dor\'enavant DGLA).
\begin{enumerate}
 \item Un ${\rm Lie}_\infty$-morphisme entre deux DGLA $ {\mathfrak g} $ et ${\mathfrak g}' $ induit une application entre \'el\'ements de Maurer-Cartan de $ {\mathfrak g} $ et \'el\'ements
 de Maurer-Cartan de $ {\mathfrak g}' $ \cite{Kont,LPV,Vallette}.
 \item Avec un \'el\'ement de Maurer-Cartan d'une DGLA $ {\mathfrak g} $, on peut d\'efinir une alg\`ebre de Leibniz diff\'erentielle gradu\'ee (ce que l'on appellera DGLeibA) 
 sur $ {\mathfrak g} [1]$, dont le crochet est ce que l'on appelle le crochet d\'eriv\'e \cite{YKS,Loday}.
\end{enumerate}
Donnons-nous un ${\rm Lie}_\infty$-morphisme $\Phi$ entre deux DGLA $ {\mathfrak g} $ 
 et ${\mathfrak g}' $ et un \'el\'ement de Maurer-Cartan $\alpha$ de $ {\mathfrak g}$. Par le premier point, on peut lui associer un \'el\'ement de Maurer-Cartan $\beta$
 de  ${\mathfrak g}' $. Par le second point, on peut utiliser $\alpha$ et $\beta$ pour construire des DGLeibA. 
 A notre connaissance, personne ne s'est pos\'e la question suivante~: existe t-il un morphisme Leibniz-infini de la DGLeibA $ {\mathfrak g} [1]$
 vers la DGLeibA ${\mathfrak g}'  [1]$~?
Pour donner un sens pr\'ecis \`a cette question, il nous faut d\'efinir ce que l'on entend par  morphisme Leibniz-infini, ce que l'on appellera un ${\rm Leib}_\infty$-morphisme.
Nous utilisons ici des id\'es similaires \`a \cite{AP}.
Rappelons d'abord qu'une DGLA $({\mathfrak g}, [\cdot,\cdot], d) $ induit sur la cog\`ebre $(S^+({\mathfrak g}[-1]), \Delta)$ une cod\'erivation $Q$ de carr\'e nul.
Celle-ci a deux coefficients de Taylor \'eventuellement non nuls: le premier coefficient de Taylor est donn\'e par $d $, et le second est le donn\'e par le crochet.
R\'eciproquement, toute cod\'erivation $Q$ de carr\'e nul sur la cog\`ebre gradu\'ee $(S^+({\mathfrak g}[-1]), \Delta)$ dont les coefficients de Taylor 
sont tous nuls sauf \'eventuellement le lin\'eaire et le quadratique induit une structure de DGLA.

Pour les alg\`ebres de Leibniz gradu\'ee, il faudra remplacer alg\`ebre sym\'etrique par alg\`ebre tensorielle $T({\mathfrak g}[-1])$. La comultiplication est alors toujours donn\'ee par 
la formule explicite que nous rappellerons dans l'\'equation (\ref{eq:comultiplication}), et qui est, comme dans la cas sym\'etrique, l'unique morphisme d'alg\`ebres gradu\'ees
$T({\mathfrak g}[-1]) \to T({\mathfrak g}[-1]) \otimes T({\mathfrak g}[-1]) $ dont la restriction \`a ${\mathfrak g}$ est donn\'ee par $ x \mapsto x \otimes 1 + 1 \otimes x$ pour tout $x \in {\mathfrak g}$.
Une difficult\'e est que cette comultiplication n'est pas colibre sur $ {\mathfrak g}$. Il reste par contre vrai que structure de DGLeibA induit une cod\'erivation de carr\'e nul de cette cog\`ebre 
\cite{LodayHomology,LodayPier,Strobl}.
La notion de coefficient de Taylor garde un sens (non-canonique) et, dans ce cas, les deux seuls coefficients \'eventuellement non-nuls 
sont la diff\'erentielle, et le crochet de Leibniz.

Nous pouvons alors d\'efinir les morphismes ${\rm Leib}_\infty$ entre deux DGLeibA comme \'etant les morphismes de cog\`ebres qui respectent ces cod\'erivations.
Notre question a maintenant un sens pr\'ecis. Nous montrons dans cet article que la r\'eponse est oui~: il existe un ${\rm Leib}_\infty$-morphisme de la DGLeibA $ {\mathfrak g} [1]$
 vers la DGLeibA $ {\mathfrak g}' [1]$. 
 Nous en d\'ecrivons les coefficients par r\'ecurrence d'une fa\c{c}on explicite.
Qui plus est, la r\'eponse est loin d'\^etre triviale~: les formules que nous obtenons sont riches en structures et font ressortir
des propri\'et\'es subtiles des alg\`ebres de Leibniz.
 
Expliquons la construction. 
Une premi\`ere complication, comme on l'a dit, est que la cog\`ebre tensorielle $T({\mathfrak g})$ (avec ${\mathfrak g}$
un espace vectoriel gradu\'e) n'est pas colibre sur ${\mathfrak g}$. Lorsque l'on traite des alg\`ebres de Leibniz, l'identit\'e
de Jacobi sugg\`ere un choix lorsque l'on fait des battages~: il faut se contenter des sommes que nous appelons \emph{respectueuses}
et que l'on note $ \stackrel{{\bullet}}{\sum}$. Pour faire comprendre le principe, rappelons que l'identit\'e de Jacobi pour les Leibniz s'\'ecrit~:
  $$ [x_1,[x_2,x_3]] - [[x_1,x_2], x_3] -  (-1)^{|x_1||x_2|} [x_2,[x_1,x_3]] =0.$$
pour tous $ x_1 , x_2, x_3 \in {\mathfrak g}$ homog\`enes.
Dans cette formule, les signes ob\'eissent aux lois de Koszul. L'ordre dans les lesquels les indices apparaissent ob\'eit 
\`a la loi suivante~: les parenth\'esages sont d'abord des battages~: c'est-\`a-dire que dans les parenth\`eses l'ordre est respect\'e.
Mais on ne conserve dans les battages que ceux qui pr\'eservent l'ordre du plus grand \'el\'ement de la parenth\`ese.
Par exemple $[[x_2,x_3],x_1]$ ne peut appara{\^i}tre car $3  > 1 $. 
D'une mani\`ere g\'en\'erale, ces sommes sont d\'efinies ainsi~:
$$   \sum^{\bullet}_{ I_1\sqcup...\sqcup I_j =[1;n]} = \sum_{ I_1\sqcup...\sqcup I_j =[1;n] \atop I_1,...,I_j\neq\emptyset, \\  p(I_1)<...< p(I_j)} $$
\noindent o\` u $ p (I) $ est le  plus grand \'el\'ement de $I$ pour toute partie $I \subset [1;n] = \{1,\dots,n\}$.


Une autre difficult\'e est un ph\'enom\`ene assez classique quand on fait de la d\'eformation~: il faut bien choisir les cobords.
C'est-\`a-dire que pour \'etendre le morphisme au degr\'e sup\'erieur, on doit chercher une certaine quantit\'e dont l'image par une diff\'erentielle
est contrainte. \'Evidemment, on peut ajouter \`a la quantit\'e cherch\'ee un cobord plus ou moins arbitraire (en tenant compte \'eventuellement des autres contraintes).
N\'eanmoins, un mauvais choix dans cet ajout peut emp\^echer de continuer la construction \`a un ordre plus \'elev\'e.
Cela signifie en particulier que l'on ne pouvait pas r\'epondre \`a la question par de simples consid\'erations cohomologiques.

Une derni\`ere subtilit\'e est qu'il faut briser la sym\'etrie. Par exemple, une partie des coefficients de Taylor sera donn\'ee par une application de l'alg\`ebre tensorielle vers l'alg\`ebre sym\'etrique
qui est d\'efinie ainsi~:
$$  x_1 \otimes  \cdots \otimes x_n \mapsto (-1)^{|x_1|+ \dots + |x_n|} Q_1 (x_1) \cdots Q_1 (x_{n-1}) \cdot x_n $$
pour tous $x_1, \dots,x_n \in {\mathfrak g}$.

D\'ecrivons maintenant le th\'eor\`eme principal de notre article. Soit ${\mathcal F}$ un  ${\rm Lie}_\infty$ morphisme d'une DGLA $ {\mathfrak g}$
vers une DGLA ${\mathfrak g}'$. 
Soit $\alpha$ un \'el\'ement de Maurer-Cartan de $ {\mathfrak g}$ et $\beta$ son image par ${\mathcal F} $.

Nous allons donner des formules tr\`es explicites car nous pensons que le lecteur expert peut comprendre la logique de notre construction en se contentant de regarder celles-ci.
Nous reprendrons les constructions plus au d\'etail tout au long du texte.
On note les d\'eriv\'ees d'ordre $n$ en $\alpha$ par $T_\alpha^n {\mathcal F}$ (voir section \ref{sec:AlgLieHomot} pour un rappel de cette notion). 
On d\'efinit une suite d'applications ${\mathcal B}_n^j~: \otimes^n {\mathfrak g} \to {\mathfrak g}' $ avec $ n \geq 1$ et $  j \geq 0 $
par r\'ecurrence. Les premiers termes sont construits ainsi~:
$$ {\mathcal B}_1^0 (x_1) = T_\alpha^1 {\mathcal F} (x_1) .$$
Ensuite, pour tout entier naturel $n$ sup\'erieur ou \'egal \` a $2$ et pour tous $x_1, \dots,x_n \in {\mathfrak g}$, on impose~:
$$    {\mathcal B }_n^0 (x_1 \otimes...\otimes x_n )= (-1)^{\vert x_1 \vert +...+ \vert x_{n-1} \vert} \ T_\alpha^n {\mathcal F} \big(Q_{\alpha,1}(x_1) \cdots Q_{\alpha,1}(x_{n-1}) \cdot x_n \big). $$
Ici, $Q_{\alpha,1} = Q_1 + Q_2(\alpha , \cdot)  $ est la diff\'erentielle associ\'e \`a l'\'el\'ement de Maurer-Cartan $\alpha$ ($Q_1$ \'etant la diff\'erentielle de ${\mathfrak g}$
et $Q_2$ son crochet de Lie gradu\'e sym\'etrique).
Enfin, on construit  pour tout $ n \geq 3$ et $ j \in \{1,...,n-2\} $ une relation de r\'ecurrence sur $j$ par~:
$$ {\mathcal B}_n^j (x_1 \otimes...\otimes x_n )= \frac{1}{j} \sum^{\bullet}_{ I\sqcup J=[1;n] } \varepsilon_x (I,J) \  \  (-1)^{\vert x_I \vert } (\vert I \vert -1) \sum_{k=0}^{j-1} 
Q'_2 \big({\mathcal B}^k_{\vert I  \vert} (x_I).  {\mathcal B}^{j-k-1}_{\vert J\vert} (x_J)\big).$$
Ici, $Q_2'$ est le crochet de Lie gradu\'e sym\'etrique de $ {\mathfrak g}'$ (\`a ne pas confondre avec le crochet d\'eriv\'e).
Le signe $ \varepsilon_x (I,J)$ ob\'eit \`a la r\`egle de Koszul.
dans la formule ci-dessus, $|I|$ est le cardinal de $I \subset \{1 , \dots,n\}$.
On convient que $ {\mathcal B}_n^j = 0 $ pour $ j \geq n-1$.

 \begin{theo*} 
  Soit ${\mathcal F}$ un  ${\rm Lie}_\infty$ morphisme d'une DGLA $ {\mathfrak g}$
vers une DGLA ${\mathfrak g}'$. Soit $\alpha $ un \'el\'ement de Maurer-Cartan 
et $ {\mathfrak g} $ et $\beta$ son image par $ {\mathcal F}$.
 
Les applications $({\mathcal B}_n)_{ n \geq 1}$ d\'efinies par $ {\mathcal B}_n := \sum_{j \geq 0 }  {\mathcal B}_n^j$  sont les coefficients de Taylor d'un 
$ {\rm Leib}_\infty$-morphisme entre les alg\`ebres de Leibniz  diff\'erentielles  gradu\'ees associ\'ees \`a $\alpha $ et $ \beta$.
 \end{theo*}

On voit bien que les formules ci-dessus brisent la sym\'etrie et qu'elles font appel \`a la notion de somme respectueuse.
Qui plus est, comme mentionn\'e ci-dessus, le choix des cobords est important.
Par exemple si on ajoute \`a $ {\mathcal B}_2$ un terme du type $ [{\mathcal B}_1 (x_1), {\mathcal B}_1 (x_2)]$,
on obtient toujours un prolongement d'ordre $2$ mais on peut montrer qu'on ne peut pas le prolonger \`a l'ordre $3$.
On voit aussi que ce  ${\rm Leib}_\infty$-morphisme n'est pas  obtenu comme une composition du ${\rm Lie}_\infty$-morphisme $ {\mathcal F}$ 
avec des d'op\'erations entre les alg\`ebres sym\'etriques et tensorielles gradu\'ees.

Comme application, nous retrouvons une formule due \`a Dominique Manchon \cite{Manchon} qui montre que la d\'eriv\'ee de l'application de formalit\'e de Kontsevich n'est pas un morphisme d'alg\`ebres 
de Lie formelles et que le d\'efaut est donn\'e par la d\'eriv\'ee seconde.

\section{Les alg\`ebres de Leibniz gradu\'ees diff\'erentielles}\label{admstr}

Commen\c{c}ons par d\'efinir l'objet principal de notre \'etude.

\begin{defi} 

\noindent 

\begin{enumerate}
\item   Une \emph{alg\`ebre de Leibniz gradu\'ee diff\'erentielle} est un triplet $({\mathfrak g}, [.,.], d) $  o\`u ${\mathfrak g} $ est un espace vectoriel gradu\'e $ {\mathfrak g}=\oplus_{i \in {\mathbb Z}} {\mathfrak g}_i$
muni 
\begin{enumerate}
\item[a)] d'un crochet $ [.,.]$ (=application bilin\'eaire $  {\mathfrak g}  \times {\mathfrak g} \mapsto {\mathfrak g}$ de degr\'e $0$)
\item[b)] d'un endomorphisme $d :  {\mathfrak g} \mapsto  {\mathfrak g}$ de degr\'e $+1$ 
\end{enumerate}
tels que et tels que pour tous $x,y,z$ homog\`enes~:
 $$ \begin{array}{rcl}[x, [y,z]] & = & [[x,y],z] + (-1)^{ \vert x\vert \vert y\vert} [y, [x,z]]  \\
      d([x,y])   & = &   [dx,y] + (-1)^{\vert x\vert} [x,dy]  \\
       d^2 & = & 0
    \end{array} $$

\item Une \emph{alg\`ebre de Lie gradu\'ee diff\'erentielle} est une alg\` ebre de Leibniz gradu\'ee diff\'erentielle  pour laquelle le crochet est antisym\'etrique gradu\'e~:

$$ [x,y] = -(-1)^{\vert x\vert \vert y\vert} [y,x].$$

\item Si ${\mathfrak g} $ est une alg\` ebre de Lie gradu\'ee diff\'erentielle. Un \'el\'ement $ \alpha $ de $ {\mathfrak g}_1 $ est dit de \emph{Maurer-Cartan} s'il v\'erifie~:

$$ d\alpha - \frac{1}{2} [ \alpha , \alpha ] = 0.$$

\end{enumerate}
\end{defi}

\noindent  Pour $ \alpha \in {\mathfrak g}_1$ , on pose
$$ d_\alpha :=[\alpha, . ] -d = {\rm ad}_\alpha -d .$$

\begin{lem}
\label{lem:MC}
Soit $({\mathfrak g}, [\cdot,\cdot],d) $ une alg\`ebre de Leibniz gradu\'ee diff\'erentielle.
Pour tout $\alpha \in {\mathfrak g} $ de degr\'e $1$,
les conditions suivantes sont \'equivalentes~:
\begin{enumerate}
 \item[(i)] $ (d_\alpha )^2 = 0 $
 \item[(ii)] pour tout $x \in {\mathfrak g}$, la relation $   [d\alpha - \frac{1}{2} [ \alpha , \alpha ], x] =0 $ est satisfaite.
\end{enumerate}
En particulier, ces conditions sont v\'erifi\'ees si $({\mathfrak g}, [\cdot,\cdot],d) $ est une alg\` ebre de Lie gradu\'ee diff\'erentielle
et $ \alpha$ est de Maurer-Cartan.
 \end{lem}

\noindent   On pose pour $x$, $y$ $ \in {\mathfrak g} $~:
$$ [x,y]_\alpha = [ (-1)^{\vert x\vert} d_\alpha (x) , y].$$

\noindent On adopte la notation habituelle sur le d\'ecalage de graduation. Ainsi, si $ V$ est un espace vectoriel gradu\'e et $n$ un entier relatif, 
un vecteur est de degr\'e $j$ dans  $ V[n] $ s'il est de degr\'e $j-n$ dans $V$.
Le crochet $[\cdot,\cdot]_\alpha$ d\'efini ci-dessus est de degr\'e $0$ sur $ { \mathfrak g} [1]$.

\begin{prop}  [Crochet d\'eriv\'e  \cite{YKS,Loday}] 
\noindent 
Soit $({\mathfrak g}, [\cdot,\cdot],d) $ une alg\` ebre de Leibniz gradu\'ee diff\'erentielle.
Pour tout $\alpha \in {\mathfrak g} $ de degr\'e $1$ qui satisfait les conditions \'equivalentes du Lemme \ref{lem:MC},
le triplet  $ ({ \mathfrak g} [1], [.,.]_\alpha , d_\alpha) $ est une alg\` ebre de Leibniz gradu\'ee diff\'erentielle.
\end{prop}

Cette proposition est en g\'en\'eral \'enonc\'ee pour  $\alpha $ un \'el\'ement Maurer-Cartan d'une alg\` ebre de Lie gradu\'ee diff\'erentielle $ ({ \mathfrak g}, [.,.], d) $, 
mais la g\'en\'eralisation ci-dessus est  \'evidente.

\section{Alg\`ebres de Lie et de Leibniz homotopiques}
\subsection{Lie$_\infty$-alg\`ebres}

\label{sec:AlgLieHomot}

On rappelle bri\`evement quelques notions d\'esormais devenues classiques sur les alg\`ebres de Lie \`a homotopie pr\`es, ou alg\`ebres $L_\infty$.
Soit $ V $ un espace vectoriel gradu\'e, l'alg\` ebre sym\'etrique de $V$ est d\'efinie par 
  $$ S (V) := T(V) / \langle x \otimes y- (-1)^{\vert x \vert \vert y\vert } y\otimes x \rangle.  $$  
  On appelle $ S^n (V) $ l'image de $ V^{\otimes n}  $ dans ce quotient
  et on d\'efinit une graduation par $ S(V) = \oplus_{n \geq 0} S^n(V)$. On consid\`ere aussi 
  $$ S^{+} (V) := \oplus_{n\geq 1} S^n (V).$$

  \noindent L'alg\`ebre  $ S^{+} (V)$  est munie d'une comultiplication coassociative~:
$$ \Delta (x_1 ... x_n) = \sum_{ I\sqcup J=[1;n]\atop I,J\neq\emptyset}\varepsilon_x (I,J) 
x_I\otimes x_J,  $$
o\`u $ \varepsilon_x (I,J) $ d\'esigne la signature de l'effet sur les $x_i$ impairs de la permutation appel\'ee
"battement" ou "battage" consistant \`a ranger d'abord les \'el\'ements de $ I$ en ordre, puis ceux de $ J$.

\begin{theo} \cite{AMM}-\cite{Kont}.
Soient $i$ un entier et $V$, $V'$ des espaces vectoriels  gradu\'es.
    Consid\'erons deux suites d'applications lin\'eaires $ Q_n ~: S^n (V) \longrightarrow V $ et 
$ {\mathcal F}_n~: S^n (V) \longrightarrow V' $ de degr\'es respectifs $i$ et $0$. Alors il existe une unique cod\'erivation $Q$ de degr\'e $i$ de $S^{+} (V) $ et un unique morphisme de cog\` ebres 
${\mathcal F}~: S^+ (V) \longrightarrow  S^+(V' )$ dont les $Q_n$ et les ${\mathcal F}_n$ sont les coefficients de Taylor respectifs. De plus, $Q$ et ${\mathcal F} $ sont donn\'es par~:
$$ Q (x_1 \dots x_n) =   \sum_{ I\sqcup J=[1;n]\atop I,J\neq\emptyset}\varepsilon_x (I,J) 
 Q_I (x_I). x_J   $$
tandis que
$$ {\mathcal F} (x_1 \dots x_n) =  \sum_{j\geq 1} \frac{1}{j!} \sum_{ I_1\sqcup...\sqcup I_j =[1;n]\atop I_1,\dots ,I_j\neq\emptyset}\varepsilon_x (I_1,\dots ,I_j) \, \, \,
 {\mathcal F}_{\vert I_1 \vert } (x_{I_1}) \dots {\mathcal F}_{\vert I_j \vert } (x_{I_j})  .$$
\end{theo}

\begin{defi}

\noindent 
\begin{enumerate}

\item  Une alg\` ebre de  Lie homotopique ou $ L_\infty$-alg\` ebre est une paire $(V,Q) $, o\`u  $V$ est un espace vectoriel gradu\'e, et o\`u $Q$
 est une cod\'erivation   de degr\'e $1$ de la cog\`ebre $  (S^+ (V), \Delta) $  v\'erifiant $ [Q,Q] =0$.
On rappelle que $[Q,Q]=2 Q^2$.

\item  Un $ L_\infty$-morphisme  entre alg\` ebres de  Lie homotopiques est un morphisme de cog\` ebres gradu\'ees:
$$ {\mathcal F}~: S^+ (V ) \rightarrow S^+(V' ) $$
 v\'erifiant~:
$$ {\mathcal F}\circ Q = Q' \circ {\mathcal F}$$

\item  La d\'eriv\'ee d'ordre $k$ d'un $ L_\infty$-morphisme ${\mathcal F}$ en un point  $ \alpha $ 
de degr\'e $0$ et v\'erifiant $ Q (e^\alpha -1) =0 $ est d\'efinie  par~:
$$ T_\alpha^k {\mathcal F} (x_1, \cdots, x_k) = \sum_{n \geq 0} \frac{1}{n!} {\mathcal F}_{n+k} (x_1  \cdots x_k. \alpha  \cdots \alpha) $$
\end{enumerate}
\end{defi}

\noindent Le cas particulier $ V = {\mathfrak g} [-1]$  o\` u $ ({\mathfrak g}, [, ], d) $ est  une alg\` ebre de Lie diff\'erentielle gradu\'ee est tr\` es important et admet des applications tr\` es int\'eressantes dans la th\'eorie de la quantification par d\'eformation.  Dans ce cas,
le champ de vecteurs $ Q$ sur $ S^{+} ({\mathfrak g} [-1] )$ a  des
coefficients de Taylor  nuls sauf les deux premiers~:

$$ Q_1 (x) = (-1)^{\vert x \vert} dx \quad \hbox {et} \quad Q_2 (x.y) = (-1)^{\vert x\vert (\vert y\vert -1)} [x,y],$$
o\` u $ \vert x\vert $ d\'esigne le degr\'e de $ x$ dans ${\mathfrak g}$, et l'\'equation $ [Q,Q] =0 $
traduit les trois relations qui d\'efinissent une alg\` ebre de Lie gradu\'ee diff\'erentielle.

\subsection{Leib$_\infty$-alg\`ebres} 

\noindent
Comme on a vu dans la section pr\'ec\'edente, les alg\`ebres $L_\infty$, en particulier les DGLA, sont encod\'ees par des cod\'erivations de carr\'e nul d'alg\`ebres sym\'etriques gradu\'ees.
Pour les alg\`ebres de Leibniz diff\'erentielles gradu\'ees, comme le crochet n'est plus sym\'etrique gradu\'e, ce codage ne peut plus se faire avec des
alg\`ebres sym\'etriques gradu\'ees. Il peut se faire, n\'eanmoins, avec des alg\`ebres tensorielles gradu\'ees.
Pr\'ecisons ce point.

\noindent Pour all\'eger les \'ecritures et \'eviter l'introduction de nouvelles notations,
dans la suite un sous-ensemble fini de $\N$ est toujours repr\'esent\'e par $ \{i_1, ... , i_k\} $ avec $ i_1<i_2<...<i_k $.

\noindent Soit $V$ un espace vectoriel gradu\'e, l'alg\` ebre tensorielle point\'ee $ T^+ (V) $ munie de  la comultiplication~:
\begin{equation}
 \label{eq:comultiplication}
 \Delta (x_1 \otimes... \otimes x_n) = \sum_{ I\sqcup J=[1;n]\atop I,J\neq\emptyset}\varepsilon_x (I,J) \ 
x_I\bigotimes x_J,  
\end{equation}
 
\noindent est une cog\` ebre cocommutative. Elle n'est pas colibre sur $ V$. 
La notion de cod\'erivation de carr\'e nul garde  \'evidemment un sens et on d\'efinit les alg\`ebres de Leibniz infinies
comme \'etant les paires $(V,Q) $, o\`u  $V$ est un espace vectoriel gradu\'e, et o\`u $Q$
 est une cod\'erivation   de degr\'e $1$ de la cog\`ebre $  (T^+ (V),\Delta) $  v\'erifiant $ [Q,Q] =0$.

 \noindent
Les morphismes entre tels objets sont les morphismes de cog\`ebres qui respectent les cod\'erivations. 

\noindent
La notion de coefficients de Taylor demande \`a \^etre pr\'ecis\'ee, car $T^+ (V) $ n'est pas colibre.
Cependant, on peut choisir une fa{\c{c}}on de construire des cod\'erivations  et des morphismes 
de cog\` ebres en partant d'une suite d'applications comme dans le cas de l'alg\` ebre sym\'etrique \cite{Kont}. 
Le choix pour les cod\'erivations est motiv\'e par la construction analogue au cas sym\'etrique d'un champ de vecteur et celui pour les morphismes nous para\^{i}t naturel. 

On a alors la proposition suivante~:

\begin{prop}
\label{prop:toutsetend}
 Soit $V$ un espace vectoriel  gradu\'e.
\begin{enumerate}
 \item 
  Consid\'erons une suite d'applications lin\'eaires $ Q_n~: T^n (V) \longrightarrow V $ de degr\'e $q \in {\mathbb Z}$. 
Consid\'erons l'application   $Q$ de $T^{+} (V) $ dans $T^{+} (V) $ d\'efinie par 

$$ Q (x_1 \otimes ...\otimes x_n) =   \sum_{k\geq 1}\sum_{ i_1<...<i_k }\varepsilon_{x,i_1,...,i_k} 
  \ x_1 \otimes ...\otimes \check{x_{i_1}} \otimes ...\otimes Q_k(x_{i_1} \otimes...\otimes  x_{i_k}) \otimes \check{x_{i_k}} \otimes ... \otimes x_n $$
\noindent  o\` u le signe ob\'eit aux r\`egles usuelles \`a condition de consid\'erer que $Q_k$ est de de degr\'e $q$, c'est-\`a-dire~:
$$ \varepsilon_{x,i_1,...,i_k}  := \varepsilon_{x'} (\{i_1,...,i_k\} , \{0,1,...,n\} \setminus \{i_1,...,i_k\} ) $$
\noindent avec $ x'_0 = Q_k $ est de degr\'e $q$ et   $ x'_l = x_l$ pour  tout $l \in \{1,...,n\} $.
L'application $Q$ est une cod\'erivation  de degr\'e $q$ dont les projections sur $V$ sont donn\'ees par les $Q_n$.

\item
Soit $ {\mathcal F}_n~: T^n (V) \longrightarrow V' $ une suite d'applications de degr\'e~$0$.
 L'application   ${\mathcal F}$ de $T^{+} (V) $ dans $T^{+} (V') $ d\'efinie par
$$ {\mathcal F} (x_1\otimes \dots \otimes x_n) =   \sum^{\bullet}_{ I_1\sqcup \dots \sqcup I_j =[1;n]}\varepsilon_x (I_1,\dots ,I_j) 
 {\mathcal F}_{\vert I_1 \vert } (x_{I_1})\otimes \dots \otimes {\mathcal F}_{\vert I_j \vert } (x_{I_j})  $$

\noindent est un morphisme de cog\` ebres dont les projections sur $V'$ sont donn\'ees par les $({\mathcal F}_n)_{n \in {\mathbb N}}$.
\end{enumerate}
\end{prop}

On qualifiera ces cod\'erivations et ces morphismes de \emph{bien faites} et on continuera \` a appeler, pour tout $ n \geq 1$, \emph{coefficients de Taylor}
la projection $Q_n : c \to V $ et $ {\mathcal F}_n : : V^{\otimes n} \to V' $ sur $V$ ou $V'$ de leur restriction \`a $ : V^{\otimes n} $.

\begin{prop}
\label{prop:bienfaites} 
Soient $V$ et $V'$ deux espaces vectoriels gradu\'es, $Q$ et $Q'$ des d\'erivations bien faites sur $V$ et $V'$ et $ {\mathcal F} : T(V) \to T(V')$ un morphisme bien fait.
Alors $ Q'  \circ{\mathcal F} = {\mathcal F}  \circ Q$ si et seulement si la projection sur $V'$ de la quantit\'e $ Q'  \circ{\mathcal F} - {\mathcal F}  \circ Q$ est nulle.
Dans le cas o\`u $Q = Q_1 +Q_2 $ et $Q'=Q_1' + Q_2'$, cela donne comme condition n\'ecessaire et suffisante les relations:
 \begin{eqnarray*}   {\mathcal F}_n \circ Q_{1}  \ \  (x_1 \otimes \cdots \otimes x_n )  + {\mathcal F}_{n-1}   \circ Q_{2}  \ \  (x_1 \otimes \cdots \otimes x_n )  &=& \\
 Q_{1}' \circ  {\mathcal F}_n   \ \ 
(x_1 \otimes \cdots \otimes x_n )  &  & \\
+ \sum^{\bullet}_{ I\sqcup J=[1;n] }  \varepsilon_x (I,J) \ \ Q'_{2}
 \left(  {\mathcal F}_{|I|} (x_I)  \otimes  {\mathcal F}_{|J|} (x_J)  \right) & &
\end{eqnarray*}
\end{prop}

\noindent Supposons maintenant que $ V = {\mathfrak g} [-1]$ o\` u $  {\mathfrak g} $ est un espace vectoriel muni d'un crochet $ [, ] $ de degr\'e $0$ et d'un endomorphisme $d$ de degr\'e $1$.
Pour $ x, y \in {\mathfrak g} $ posons~: 

$$ Q_1 (x) = (-1)^{\vert x \vert} dx \quad \hbox {et} \quad Q_2 (x\otimes y) = (-1)^{\vert x\vert (\vert y\vert -1)} [x,y],$$
o\` u $ \vert x\vert $ d\'esigne le degr\'e de $ x$ dans ${\mathfrak g}$.

\noindent Appelons $ Q$ la cod\'erivation construite \` a partir de $ Q_1$ et $Q_2$ comme dans la proposition pr\'ec\'edente. On a bien 
$$ Q^2 (x) = Q_1^2 (x) $$
\noindent  et 
$$ Q^2 (x\otimes y) = Q_1^2 (x) \otimes y + x \otimes Q_1^2 (y) + Q_1 (Q_2 (x\otimes y) + Q_2 (Q_1(x) \otimes y) + (-1)^{\vert x \vert -1} Q_2 ( x\otimes Q_1(y)),$$

\noindent Un calcul direct montre que pour tout  entier naturel $n$ sup\'erieur \` a $3$ et pour tout $x_1,..., x_n$ dans ${\mathfrak g} [-1] $ on a~:

\begin{flalign*} & Q^2 (x_1\otimes ... \otimes x_n) = \sum_{i=1}^n x_1 \otimes ...\otimes Q_1^2 (x_i)\otimes ...\otimes x_n + \sum_{ i<j }\varepsilon_{x} (\{i,j\},\{1,...,n\}  \setminus \{i,j\} ) \\
 &  x_1 \otimes ...\otimes \check{x_{i}} \otimes ...\otimes \Big(Q_1 \big(Q_2(x_i \otimes x_j)\big) + Q_2\big(Q_1(x_i)\otimes x_j\big) + (-1)^{\vert x_i\vert -1} Q_2\big(x_i \otimes Q_1(x_j)\big)\Big)  \otimes  \check{x_{j}} \otimes ... \otimes x_n \\
 & + \sum_{ i<j<k }\varepsilon_{x} (\{i,j,k\},\{1,...,n\}  \setminus \{i,j,k\} )\  x_1 \otimes ...\otimes \check{x_{i}} \otimes ...\otimes \check{x_{j}} \otimes...  \\
 & \otimes  \Big(Q_2 \big(Q_2(x_i \otimes x_j)\otimes x_k\big) + (-1)^{ \vert x_i \vert -1 }Q_2\big(x_i\otimes Q_2(x_j\otimes x_k)\big) + (-1)^{\vert x_i \vert (\vert x_j\vert -1)} Q_2\big(x_j \otimes Q_2(x_i \otimes x_k)\big)\Big) \\
 &  \otimes  \check{x_{k}} \otimes ... \otimes x_n.
 \end{flalign*}

\noindent Pour voir que les termes du type~:

$$ x_1 \otimes ... \otimes Q_1(x_i) \otimes \check{x_{i}} \otimes ...\otimes \check{x_{j}} \otimes ...\otimes Q_2(x_j\otimes x_k) \otimes \check{x_{k}}\otimes ...\otimes x_n   $$

\noindent se simplifient, il suffit de remarquer qu'un tel terme appara\^it quand on applique $Q_1$ \` a $ Q_2 (x_1 \otimes ... \otimes x_n)$ sans que $Q_1$ et $Q_2$ ne se 
rencontrent mais aussi quand on applique $Q_2$ \` a $ Q_1 (x_1 \otimes ... \otimes x_n)$ sauf que  dans ce cas $Q_2$ est pass\'e au dessus de $Q_1$, ce qui fait na\^itre un signe oppos\'e.

\noindent On d\'eduit alors la proposition suivante, qui g\'en\'eralise quelques r\'esultats de \cite{LodayHomology,LodayPier}~:

\begin{prop} 
Soient $ ({\mathfrak g} , [\cdot,\cdot], d) $, un espace vectoriel gradu\'e, muni d'une application bilin\'eaire de degr\'e $0$ et une application lin\'eaire de degr\'e $+1$.
Soit $Q$ comme ci-dessus.
Alors $ Q^2 = 0 $ si et seulement si $ ({\mathfrak g} , [\cdot,\cdot], d) $ est une alg\` ebre de Leibniz diff\'erentielle gradu\'ee.
\end{prop}

\section{De Lie \` a Leibniz}

\noindent  Soient $ ({\mathfrak g}, [, ], d) $ et $ ({\mathfrak g}', [, ]', d') $ deux 
alg\` ebres de 
Lie gradu\'ees diff\'erentielles .On a donc deux $ L_\infty$-alg\` ebres  $ (S^{+} ({\mathfrak g} [-1] ), Q)$ et $ (S^{+} ({\mathfrak g'} [-1] ), Q')$. Supposons qu'il existe un $ L_\infty$-morphisme 
${\mathcal F}$ entre ces deux $ L_\infty$-alg\` ebres.

\noindent  Soit $ \alpha $ un \'el\'ement de Maurer-Cartan de $ {\mathfrak g}$, on a donc~:

$$ Q_1 (\alpha) + \frac{1}{2} Q_2 (\alpha . \alpha) = 0.$$

\noindent  Posons~:

$$ \beta = \sum_{n=1}^{+\infty} \frac{1}{n!}{\mathcal F}_n ( \alpha ... \alpha ) = pr ( {\mathcal F} (e^\alpha -1)) \ \in {\mathfrak g}'_1 .$$

\noindent  Les relations suivantes d\'ecoulent du lemme \ref{lem:ManchonKont} rappel\'e ci-dessous~:

$$ {\mathcal F} ( e^\alpha -1) = e^\beta -1 .$$

$$ (Q'_1 (\beta ) + \frac{1}{2} Q'_2 (\beta . \beta ) ) e^\beta = Q' ( e^\beta -1) = Q' ( {\mathcal F} (e^\alpha -1)) = {\mathcal F } ( Q (e^\alpha -1)  ) = 0. $$

\noindent Il en suit que $ \beta $ est un \'el\'ement de Maurer-Cartan de ${\mathfrak g}'$.

\noindent  En consid\'erant les structures d\'eriv\'ees correspondantes \` a $\alpha $ et $\beta $, on obtient  deux   alg\` ebres de 
Leibniz gradu\'ees diff\'erentielles $ ({ \mathfrak g} [1], [.,.]_\alpha , d_\alpha) $ et $ ({ \mathfrak g}' [1], [.,.]'_\beta , d'_\beta) $ (et par suite deux alg\` ebres
de Leibniz 'homotopiques')~: 

$$ (T^+ ( {\mathfrak g}[1][-1] , Q_\alpha ) = (T^+ ({\mathfrak g}) , Q_\alpha )\quad \hbox {et} \quad(T^+ ({\mathfrak g}') , Q_\beta ).$$

\noindent Pr\'ecis\'ement, on a~: 

$$ Q_{\alpha, 1} ( x) = (-1)^ {\vert x\vert +1 } d_\alpha (x)  \quad , \quad Q_{\alpha,2} ( x \otimes y ) = (-1)^{ (\vert x \vert +1) \vert y \vert } [x,y]_\alpha  $$

\noindent et

$$ Q'_{\beta, 1} ( x) = (-1)^ {\vert x\vert +1 } d'_\beta (x)  \quad , \quad Q'_{\beta,2} ( x \otimes y ) = (-1)^{ (\vert x \vert +1) \vert y \vert } [x,y]_\beta' . $$

\noindent Le lemme suivant est une r\'e\'ecriture de celle donn\'ee dans \cite{Kont} et \cite{Manchon} qui utilise les notations du paragraphe 
pr\'ec\'edent et le fait que pour tout $x$ dans  $ {\mathfrak g}$
$$ Q_{\alpha,1} (x) = Q_1(x) + Q_2 (x.\alpha).$$
Le calcul se fait dans l'alg\`ebre sym\'etrique d\'ecal\'ee, c'est-\`a-dire dans $S^{+} ({\mathfrak g} [-1]) $ et $ S^{+} ({\mathfrak g}' [-1])$.
En particulier, les signes $ \varepsilon_x (I_1,...,I_j)$ sont ceux qui apparaissent dans cette alg\`ebre. 
Par exemple, $ \varepsilon_{x,2} = (-1)^{\vert x_1 \vert -1}$ si $ \vert x_1 \vert $ est le degr\'e de $ x_1$ dans $ {\mathfrak g}$.


\begin{lem} \cite{Kont}-\cite{Manchon}
\label{lem:ManchonKont}
\noindent 
Soit ${\mathcal F}: S^+( {\mathfrak g}[1]) \to  S^+( {\mathfrak g}'[1]) $ un $L_\infty$-morphisme d'une DGLA $({\mathfrak g}, [\cdot,\cdot],d) $ vers une DGLA $({\mathfrak g}', [\cdot,\cdot]',d') $.
Soit $\alpha \in {\mathfrak g} $ un \'el\'ement de Maurer-Cartan. Pour tout 
$  x_1,..., x_n \in {\mathfrak g} $,
on a~:
$$ {\mathcal F} (x_1  ...  x_n. e^\alpha ) = \Bigl( \sum_{j\geq 1}  \sum^{\bullet}_{ I_1\sqcup...\sqcup I_j =[1;n]}\varepsilon_x (I_1,...,I_j) 
\ T_\alpha^{\vert I_1 \vert } {\mathcal F} (x_{I_1})  ...   T_\alpha^{\vert I_j \vert } {\mathcal F} (x_{I_j})\Bigr) e^\beta.  $$
 \begin{flalign*}  Q (x_1 ...  x_n e^\alpha )& = \Bigl(\sum_{i=1}^n \varepsilon_{x,i} \ x_1  ... \check{x_i}  Q_{\alpha,1} (x_i)  ... x_n  \\
&+ \sum_{ i<j} \varepsilon_{x,i,j} \ x_1   ...  \check{x_i}   ... Q_2(x_i . x_j)  \check{x_j}  ...  x_n \Bigr)  e^\alpha .
\end{flalign*} 
\end{lem}

D\'efinissons maintenant une suite d'applications de $T^n({\mathfrak g}) $ dans ${\mathfrak g}'$.
 \noindent On pose
$$ {\mathcal B}_1^0 (x_1) = T_\alpha^1 {\mathcal F} (x_1) $$
\noindent et pour tout entier naturel $n$ sup\'erieur ou \'egal \` a $2$
$$    {\mathcal B }_n^0 (x_1 \otimes...\otimes x_n )= (-1)^{\vert x_1 \vert +...+ \vert x_{n-1} \vert} \ T_\alpha^n {\mathcal F} \big(Q_{\alpha,1}(x_1)...Q_{\alpha,1}(x_{n-1}).x_n\big). $$
\noindent C'est-\` a-dire, en utilisant les notations pr\'ec\'edentes
$$    {\mathcal B }_n^0 (x_1 \otimes...\otimes x_n )=  \varepsilon_{x,n} \ T_\alpha^n {\mathcal F} \big(Q_{\alpha,1}(x_1)...Q_{\alpha,1}(x_{n-1}).x_n\big). $$ 

\noindent Pour tout $ n \geq 3$,  on d\'efinit ${\mathcal B}_n^j $ par une relation de r\'ecurrence sur $ j \in \{1,...,n-2\} $ de la fa\c{c}on suivante~:
$$ {\mathcal B}_n^j (x_1 \otimes...\otimes x_n )= \frac{1}{j} \sum^{\bullet}_{ I\sqcup J=[1;n] } \varepsilon_x (I,J) \  \  (-1)^{\vert x_I \vert } (\vert I \vert -1) \sum_{k=0}^{j-1} 
Q'_2 \big({\mathcal B}^k_{\vert I  \vert} (x_I).  {\mathcal B}^{j-k-1}_{\vert J\vert} (x_J)\big).$$
Ici, $|I|$ d\'esigne le cardinal de $I \subset \{1, \dots,n\}$.

\noindent On convient de poser $ {\mathcal B}_n^j = 0 $ dans les autres cas.
Finalement, on pose 

$$ {\mathcal B}_n = \sum_{ j \geq 0} {\mathcal B}_n^j.$$

\noindent
Voici le r\'esultat principal de notre article~:

 \begin{theo}
 \label{theo:principal}
 Soit ${\mathcal F } $ un $L_\infty$-morphisme d'une DGLA $ ({\mathfrak g} ,  [ \cdot,\cdot] , d )$ vers une DGLA $ ({\mathfrak g} ',  [ \cdot,\cdot]' , d' )$.
Soit $\alpha \in {\mathfrak g}_1$ un \'el\'ement de Maurer-Cartan de  $ ({\mathfrak g} ,  [ \cdot,\cdot] , d )$  et $\beta  \in {\mathfrak g} '_1$ l'\'el\'ement de Maurer-Cartan  
image dans~$ {\mathfrak g}'$.
 
 \noindent La suite d'applications $({\mathcal B}_n)_{ n \geq 1} $ est la suite des   coefficients de Taylor d'un 
${\rm Leib}_\infty$-morphisme $ {\mathcal B}:T^+ ({\mathfrak g}) \to T^+ ({\mathfrak g}') $ entre les alg\`ebres de Leibniz diff\'erentielles gradu\'ees 
$ ({\mathfrak g} ,  [ \cdot,\cdot]_\alpha , d_\alpha )$ et $ ({\mathfrak g}',  [ \cdot,\cdot]_\beta' , d_\beta' )$.
 \end{theo}
 
 \noindent La preuve de ce th\'eor\` eme est une cons\'equence des deux propositions \ref{prop:Cn0} et \ref{prop:Cnj} suivantes. Consid\'erons la quantit\'e:
 $$
  {\mathcal C}_n^0 ( x_1\otimes ... \otimes x_n)  = \big({\mathcal B }_n^0 \circ Q_{\alpha, 1} + {\mathcal B }_{n-1}^0 \circ Q_{\alpha,2} - Q'_{\beta,1} \circ {\mathcal B}_n^0 \big) (x_1\otimes ... \otimes x_n) \\
 $$

\begin{prop}
\label{prop:Cn0}
 Pour tout entier naturel $n$ sup\'erieur ou \'egal \` a $2$ et pour tout $x_1,...,x_n $ dans ${\mathfrak g}$
  \begin{flalign*} {\mathcal C}_n^0 ( x_1\otimes ... \otimes x_n)  = & 
  \sum^{\bullet}_{ I\sqcup J=[1;n] \atop \vert I \vert = 1 }  \varepsilon_x (I,J) Q'_{\beta,2} \big({\mathcal B}_1^0 (x_I) \otimes {\mathcal B}_{n-1}^0 (x_J) \big) \\ 
 + & \sum^{\bullet}_{ I\sqcup J=[1;n] \atop \vert I \vert \geq 2 } \ \ (-1)^{\vert x_I \vert } \varepsilon_x (I,J) Q'_2 
 \Big({\mathcal B}_{\vert I\vert}^0 \big( Q_{\alpha, 1} (x_I)\big) \otimes {\mathcal B}_{\vert J\vert}^0 (x_J)\Big)   \end{flalign*}  
\end{prop}

\begin{preuve}
\noindent En utilisant le lemme \ref{lem:ManchonKont} et en remarquant que $ Q_{\alpha,1}^2 = 0$, on peut voir que la projection de la relation (\'ecrite dans $S^{+} ({\mathfrak g'} [-1] )$):
$$ {\mathcal F} \circ Q \big( Q_{\alpha,1} (x_1)...Q_{\alpha,1} (x_{n-1}).x_n e^\alpha \big) = Q' \circ {\mathcal F} \big( Q_{\alpha,1} (x_1)...Q_{\alpha,1} (x_{n-1}).x_n e^\alpha \big)$$
\noindent se traduit par~:

\begin{flalign*}
&\varepsilon_{x',n}\  T_\alpha^n {\mathcal F} \big(Q_{\alpha,1}(x_1)...Q_{\alpha,1}(x_{n-1}).Q_{\alpha,1}(x_n)\big)\\ 
   +  & \sum_{ i<j<n} \varepsilon_{x',i,j} \ T_\alpha^{n-1} {\mathcal F} \big(Q_{\alpha,1}(x_1)...
  \check {Q_{\alpha,1}(x_i)}... Q_2\big(Q_{\alpha,1}(x_i).Q_{\alpha,1}(x_j)\big).\check{Q_{\alpha,1}(x_j)}...x_n\big)\\
  + &  \sum_{i=1}^{n-1} \varepsilon_{x',i,n} \ T_\alpha^{n-1} {\mathcal F} \Big(Q_{\alpha,1}(x_1)...
  \check {Q_{\alpha,1}(x_i)}... Q_2\big(Q_{\alpha,1}(x_i).x_n\big)\Big)\\
   = & Q'_{\beta,1} \Big(T_\alpha^n {\mathcal F} \big(Q_{\alpha,1}(x_1)...Q_{\alpha,1}(x_{n-1}).x_n\big)\Big) + 
  \sum^{\bullet}_{ I\sqcup J =[1;n]}\varepsilon_{x'} (I,J) 
Q'_2 \big(\ T_\alpha^{\vert I \vert } {\mathcal F} (x'_{I}) . T_\alpha^{\vert J \vert } {\mathcal F} (x'_{J})\big) 
\end{flalign*}

\noindent o\` u $ x'_i = Q_{\alpha,1} (x_i) $ si $ i< n$ et $ x'_n = x_n $.

\noindent Or $ Q_2\big(Q_{\alpha,1} (x_i) . Q_{\alpha,1} (x_j) \big) = - Q_{\alpha,1} \big(Q_{\alpha,2} (x_i\otimes x_j )\big)$ et $ Q_2\big(Q_{\alpha,1} (x_i) . x_n \big) = (-1)^{\vert x_i \vert}  Q_{\alpha,2} (x_i\otimes x_n )$.

\noindent En multipliant   par $ \varepsilon_{x',n}$ l'\'egalit\'e devient~:

\begin{flalign*}
&   T_\alpha^n {\mathcal F} \big(Q_{\alpha,1}(x_1)...Q_{\alpha,1}(x_{n-1}).Q_{\alpha,1}(x_n)\big)\\ 
  + &  \sum_{ i<j<n} \varepsilon_{x',i,j} (-\varepsilon_{x',n}) \ T_\alpha^{n-1} {\mathcal F} \big(Q_{\alpha,1}(x_1)...
  \check {Q_{\alpha,1}(x_i)}... Q_{\alpha,1}\big(Q_{\alpha,2}(x_i \otimes x_j)\big).\check{Q_{\alpha,1}(x_j)}...x_n \big)\\
  + &  \sum_{i=1}^{n-1} \varepsilon_{x',i,n} \varepsilon_{x',n} (-1)^{\vert x_i \vert }\ T_\alpha^{n-1} {\mathcal F} \big(Q_{\alpha,1}(x_1)...
  \check {Q_{\alpha,1}(x_i)}...  Q_{\alpha,2}(x_i\otimes x_n)\big) \\
   - & Q'_{\beta,1} \Big(\varepsilon_{x',n} \ T_\alpha^n {\mathcal F} \big(Q_{\alpha,1}(x_1)...Q_{\alpha,1}(x_{n-1}).x_n\big)\Big)  \\
   = & \varepsilon_{x',n} \sum^{\bullet}_{ I\sqcup J=[1;n] \atop \vert I \vert = 1 }\varepsilon_{x'} (I,J) 
Q'_2 \big(\ T_\alpha^1 {\mathcal F} (x'_I) . T_\alpha^{\vert J \vert } {\mathcal F} (x'_J)\big) + \varepsilon_{x',n} \sum^{\bullet}_{ I\sqcup J=[1;n] \atop \vert I \vert > 2 }\varepsilon_{x'} (I,J) 
Q'_2 \big(\ T_\alpha^{\vert I \vert } {\mathcal F} (x'_I) . T_\alpha^{\vert J \vert } {\mathcal F} (x'_J)\big) 
\end{flalign*}

\noindent Pour pouvoir conclure, il faut remarquer que pour tout $ i \in \{1,...,n-1\} $, ${\mathcal B}_n^0 \big(x_1\otimes ..\otimes Q_{\alpha,1}(x_i) \otimes...\otimes x_n\big) = 0$ et qu'en passant  de $S^{+} ({\mathfrak g} [-1])$ \` a $T^{+} ({\mathfrak g}) $  on a 

$$ \varepsilon_{x',n} = \varepsilon_{x,n} = (-1)^{\vert x_1 \vert +...+ \vert x_{n-1} \vert},\quad \varepsilon_{x',i,j} =\varepsilon_{x ,i,j}  $$

\noindent et puisqu'il s'agit d'une somme respectueuse 

$$ \varepsilon_{x'} (I,J) = \varepsilon_{x} (I,J). $$ 

\noindent Enfin, il faut aussi se rappeler que $  T_\alpha^1$ est un morphisme de complexes~:
$$ T_\alpha^1 \circ Q_{\alpha,1} = Q'_{\beta,1} \circ T_\alpha^1,$$
ce qui donne le r\'esultat.
\end{preuve}

 Pour tout entier naturel $j$ sup\'erieur ou \'egal \` a $1$ et pour tout $x_1,...,x_n $ ($ n \geq  j+2$) dans ${\mathfrak g}$, on introduit~:
$$
 {\mathcal C}_n^j ( x_1\otimes ... \otimes x_n)  = \big({\mathcal B }_n^j \circ Q_{\alpha, 1} + {\mathcal B }_{n-1}^j \circ Q_{\alpha,2} - Q'_{\beta,1} \circ {\mathcal B}_n^j  \big) (x_1\otimes 
 \cdots \otimes x_n) 
 $$ 

\begin{prop}
\label{prop:Cnj}
 Pour tout entier naturel $j$ sup\'erieur ou \'egal \` a $1$ et pour tous $x_1,...,x_n $ \emph{(}avec $ n \geq  j+2$\emph{)} dans ${\mathfrak g}$~:
  \begin{flalign*} \label{eq:propCnj} & {\mathcal C}_n^j ( x_1\otimes ... \otimes x_n) \\ 
 = & \sum^{\bullet}_{ I\sqcup J=[1;n] \atop \vert I \vert \geq 2 }  \varepsilon_x (I,J)
 \sum_{k=0}^{j-1} Q'_{\beta,2} \big({\mathcal B}_{\vert I \vert }^k (x_I) \otimes {\mathcal B}_{\vert J \vert }^{j-k-1} (x_J) \big) \\ 
 +  & \sum^{\bullet}_{ I\sqcup J=[1;n] \atop \vert I \vert = 1 }  \varepsilon_x (I,J)
 \sum_{k=0}^j  Q'_{\beta,2} \big({\mathcal B}_1^k (x_I) \otimes {\mathcal B}_{n-1}^{j-k} (x_J) \big) \\
 +  & \ \ \ R_n^j (x_1 \otimes... \otimes x_n) - R_n^{j+1} (x_1 \otimes... \otimes x_n) \end{flalign*} 
 
 \noindent o\`u pour tout $ m \in \mathbb{N}^* $~:
 $$ R_n^m  (x_1 \otimes... \otimes x_n) = \sum^{\bullet}_{ I\sqcup J=[1;n] \atop \vert I \vert \geq 2 }  \varepsilon_x (I,J) (-1)^{\vert x_I \vert + 1} \sum_{k=0}^{m-1} Q'_2 \big({\mathcal B}_{\vert I \vert }^k \big( Q_{\alpha , 1}(x_I)\big) .
 {\mathcal B}_{\vert J \vert }^{m-k-1} (x_J) \big) .$$
\end{prop}

\begin{rem}
\label{rmque:Cn0}
D'apr\` es la proposition \ref{prop:Cn0}, la relation suivante est satisfaite~: 
$$ {\mathcal C}_n^0 ( x_1\otimes ... \otimes x_n)  = 
  \sum^{\bullet}_{ I\sqcup J=[1;n] \atop \vert I \vert = 1 }  \varepsilon_x (I,J) Q'_{\beta,2} \big({\mathcal B}_1^0 (x_I) \otimes {\mathcal B}_{n-1}^0 (x_J) \big) -R_n^1. $$
\end{rem}

\noindent La preuve de la proposition \ref{prop:Cnj} utilise les deux lemmes suivants~:

\begin{lem}\label{lem:2rel}
Les deux relations suivantes sont satisfaites~:
\begin{flalign*} {\mathcal B}_n^j \big(Q_{\alpha,1}(x_1 \otimes...\otimes x_n )\big)& = 
  \frac{1}{j} \sum^{\bullet}_{ I\sqcup J=[1;n] } \varepsilon_x (I,J) \  \  (-1)^{\vert x_I \vert +1 } (\vert I \vert -1) \sum_{k=0}^{j-1} 
Q'_2 \big({\mathcal B}^k_{\vert I  \vert} \big(Q_{\alpha,1}(x_I)\big) . {\mathcal B}^{j-k-1}_{\vert J \vert} (x_J)\big) \\ 
& +
\frac{1}{j} \sum^{\bullet}_{ I\sqcup J=[1;n] } \varepsilon_x (I,J) \  \   (\vert I \vert -1) \sum_{k=0}^{j-1} 
Q'_2 \Big({\mathcal B}^k_{\vert I  \vert} (x_I) .{\mathcal B}^{j-k-1}_{\vert J \vert } \big(Q_{\alpha,1}(x_J)\big)\Big) \end{flalign*}
\noindent 
et~:
\begin{flalign*} {\mathcal B}_{n-1}^j \big(Q_{\alpha,2}(x_1 \otimes...\otimes x_n )\big)& = 
  \frac{1}{j} \sum^{\bullet}_{ I\sqcup J=[1;n] \atop \vert I \vert \geq 2 } \varepsilon_x (I,J) \  \  (-1)^{\vert x_I \vert +1 } (\vert I \vert -2) \sum_{k=0}^{j-1} 
Q'_2 \big({\mathcal B}^k_{\vert I  \vert-1} \big(Q_{\alpha,2}(x_I)\big) . {\mathcal B}^{j-k-1}_{\vert J \vert } (x_J)\big) \\ 
& +
\frac{1}{j} \sum^{\bullet}_{ I\sqcup J=[1;n] \atop \vert J \vert \geq 2} \varepsilon_x (I,J) \  \   (\vert I \vert -1) \sum_{k=0}^{j-1} 
Q'_2 \Big({\mathcal B}^k_{\vert I  \vert} (x_I) .{\mathcal B}^{j-k-1}_{\vert J \vert -1} \big(Q_{\alpha,2}(x_J)\big)\Big) \end{flalign*}
\end{lem}

\begin{preuve} (du lemme \ref{lem:2rel}). En utilisant la d\'efinition de $ {\mathcal B}_n^j $, on obtient~:
\begin{flalign*} {\mathcal B}_n^j \big(Q_{\alpha,1}(x_1 \otimes...\otimes x_n )\big)& = \sum_{i=1}^n \varepsilon_{x,i}  {\mathcal B}_n^j (x_1 \otimes...\otimes Q_{\alpha,1}(x_i)\otimes ...\otimes x_n )
\\& = \sum_{i=1}^n \varepsilon_{x,i}  \frac{1}{j} \sum^{\bullet}_{ I\sqcup J=[1;n] } \varepsilon_{x'} (I,J) \  \  (-1)^{\vert x'_I \vert } (\vert I \vert -1) \sum_{k=0}^{j-1} 
Q'_2 \big({\mathcal B}^k_{\vert I  \vert} (x'_I).  {\mathcal B}^{j-k-1}_{\vert J\vert} (x'_J)\big)
\\& =   \sum^{\bullet}_{ I\sqcup J=[1;n] \atop  i \in I } \varepsilon_{x,i}  \frac{1}{j} \varepsilon_{x'} (I,J) \  \  (-1)^{\vert x'_I \vert } (\vert I \vert -1) \sum_{k=0}^{j-1} 
Q'_2 \big({\mathcal B}^k_{\vert I  \vert} (x'_I).  {\mathcal B}^{j-k-1}_{\vert J\vert} (x'_J)\big)
\\& +  \sum^{\bullet}_{ I\sqcup J=[1;n] \atop   i \in J }  \varepsilon_{x,i}  \frac{1}{j} \varepsilon_{x'} (I,J) \  \  (-1)^{\vert x'_I \vert } (\vert I \vert -1) \sum_{k=0}^{j-1} 
Q'_2 \big({\mathcal B}^k_{\vert I  \vert} (x'_I).  {\mathcal B}^{j-k-1}_{\vert J\vert} (x'_J)\big)
\end{flalign*}

\noindent o\` u $ x'_k= x_k $ si $ k \neq i$ et $ x'_i = Q_{\alpha,1}(x_i)$. On remarque que les termes de la somme
$  \sum_{k=0}^{j-1}  Q'_2 \big({\mathcal B}^k_{\vert I  \vert} (x'_I).  {\mathcal B}^{j-k-1}_{\vert J\vert} (x'_J)\big)$, sont, au signe pr\`es, les m\^emes termes que ceux de la somme
$ \sum_{k=0}^{j-1} 
Q'_2 \big({\mathcal B}^k_{\vert I  \vert} \big(Q_{\alpha,1}(x_I)\big) . {\mathcal B}^{j-k-1}_{\vert J \vert} (x_J)\big) $.
Pour conclure, il faut v\'erifier les signes. Si $ i \in I$ alors 
$$ \varepsilon_{x'}( I,J) = (-1)^{ \sum_{ k<i, k \in J} \vert x_k \vert }\  \varepsilon_x (I,J) \quad , \quad  \vert x'_I \vert = \vert x_I \vert +1 $$
\noindent  par suite
$$ \varepsilon_{x,i} \varepsilon_{x'}( I,J) (-1)^{\vert x'_I \vert} = (-1)^{ \sum_{ k<i, k \in I} \vert x_k \vert } (-1)^{\vert x_I \vert +1} \varepsilon_x (I,J),$$
\noindent et si $ i \in J$ alors
$$ \varepsilon_{x'}( I,J) = (-1)^{ \sum_{ k>i, k \in I} \vert x_k \vert }\  \varepsilon_x (I,J) \quad , \quad  \vert x'_I \vert = \vert x_I \vert   $$
\noindent  par suite
$$ \varepsilon_{x,i} \varepsilon_{x'}( I,J) (-1)^{\vert x'_I \vert} = (-1)^{ \sum_{ k<i, k \in J} \vert x_k \vert } \varepsilon_x (I,J).$$
Ceci montre la premi\`ere formule.

\noindent Pour la seconde formule, on a aussi~:

\begin{flalign*}  &  {\mathcal B}_{n-1}^j (Q_{\alpha,2}(x_1 \otimes...\otimes x_n )) \\ & = \sum_{l<s}  \varepsilon_{x,l,s}  \, {\mathcal B}_{n-1}^j (x'_1 \otimes...\otimes x'_{n-1} )\\& = 
\sum_{l<s}  \varepsilon_{x,l,s}  \frac{1}{j} \sum^{\bullet}_{ I'\sqcup J'=[1;n-1] } \varepsilon_{x'} (I',J') \  \  (-1)^{\vert x'_{I'} \vert } (\vert I' \vert -1) \sum_{k=0}^{j-1} 
Q'_2 \big({\mathcal B}^k_{\vert I'  \vert} (x'_I).  {\mathcal B}^{j-k-1}_{\vert J'\vert} (x'_{J'})\big)
\\& = \sum_{l<s}  \varepsilon_{x,l,s}  \frac{1}{j} \sum^{\bullet}_{ I'\sqcup J'=[1;n-1 ] \atop s-1 \in I' } \varepsilon_{x'} (I',J') \  \  (-1)^{\vert x'_{I'} \vert } (\vert I' \vert -1) \sum_{k=0}^{j-1} 
Q'_2 \big({\mathcal B}^k_{\vert I'  \vert} (x'_{I'}).  {\mathcal B}^{j-k-1}_{\vert J'\vert} (x'_{J'})\big)
\\& + \sum_{l<s}  \varepsilon_{x,l,s}  \frac{1}{j} \sum^{\bullet}_{ I'\sqcup J'=[1;n-1] \atop s-1 \in J' } \varepsilon_{x'} (I',J') \  \  (-1)^{\vert x'_{I'} \vert } (\vert I' \vert -1) \sum_{k=0}^{j-1} 
Q'_2 \big({\mathcal B}^k_{\vert I'  \vert} (x'_{I'}).  {\mathcal B}^{j-k-1}_{\vert J'\vert} (x'_{J'})\big)
\end{flalign*}

\noindent o\` u 

$$ x'_1= x_1~; ...; x'_{l-1}= x_{l-1}; x'_l= x_{l+1};...x'_{s-2}= x_{s-1};   x'_{s-1} = Q_{\alpha,2}(x_l \otimes x_s); x'_s = x_{s+1};...; x'_{n-1} = x'_n.$$

\noindent Pour $ l$ et $s$ fix\'es tels que $ l<s$, \` a un sous ensemble $I'$ de $ \{ 1,..., n-1\}$ on associe un sous-ensemble $I$ de $ \{ 1,..., n \}$ de la fa\c con  suivante~:
\begin{enumerate}
 \item[a)] si $ s-1 \in I'$, alors $I$ contiendra les \'el\'ements de $I'$ inf\'erieurs strictement \` a $l$, $l$ et $s$ et les autres \'el\'ements de $I$ sont de la forme $j+1$ avec $ j \in I'$ et $ j \leq l$.
 \item[b)] sinon $I$ sera form\'e des \'el\'ements de $I'$ strictement inf\'erieurs \` a $l$ et des \'el\'ements de la forme $j+1$ avec $ j \in I'$ et $ j \leq l$.
\end{enumerate}

\noindent Il est clair que si $ s-1 \in I'$, $ \vert I \vert  = \vert I' \vert +1 $ et $ (-1)^{ \vert x'_{I'} \vert } = (-1)^{\vert x_I \vert +1} $, sinon  $ \vert I \vert  = \vert I' \vert  $ et $ (-1)^{ \vert x'_{I'} \vert } = (-1)^{\vert x_I \vert } $.

\noindent Les termes du second membre de l'\'egalit\'e pr\'ec\'edente sont alors de m\^eme nature que ceux de la formule donn\'ee dans le lemme. Pour les signes il suffit de voir que~:

\noindent  Si $ s-1 \in I'$, alors

$$ \varepsilon_{x'}( I,J) = (-1)^{ \sum_{ k<l, k \in J} \vert x_k \vert + \sum_{ l\leq k<s-1, k \in J'} \vert x_{k+1} \vert (\vert x_l \vert +1)}\  \varepsilon_x (I,J)  $$

\noindent sinon 

$$ \varepsilon_{x'}( I,J) = (-1)^{ \sum_{ k>s-1, k \in I'} \vert x_{k+1} \vert + \sum_{ l\leq k<s-1, k \in I'} \vert x_{k+1} \vert  \vert x_l \vert }\  \varepsilon_x (I,J)  $$
\end{preuve}

\begin{lem} 
\label{lem:sommes}
Les deux relations suivantes sont satisfaites~:
\begin{flalign*}  
 & \sum_{l=0}^m \Bigg( \sum^{\bullet}_{ I\sqcup J=[1;n] \atop \vert I\vert \geq 3} \varepsilon_x (I,J) \    (-1)^{\vert x_I \vert +1 } (\vert I \vert -2) \sum^{\bullet}_{I_1 \sqcup I_2 = I \atop \vert I_1 \vert \geq 2} \varepsilon_x (I_1,I_2) (-1)^{\vert x_{I_1} \vert  }
Q'_2 \bigg( Q'_2 \Big({\mathcal B}^l_{\vert I_1  \vert} \big(Q_{\alpha,1}(x_{I_1}\big)\Big) . {\mathcal B}^{m-l}_{\vert I_2 \vert} (x_{I_2})). {\mathcal B}^k_{\vert J \vert} (x_{J})\bigg) \\ 
& + \sum^{\bullet}_{ I\sqcup J=[1;n] \atop \vert J\vert \geq 2} \varepsilon_x (I,J) \      (\vert I \vert -1) \sum^{\bullet}_{J_1 \sqcup J_2 = J \atop \vert J_1 \vert \geq 2} \varepsilon_x (J_1,J_2) (-1)^{\vert x_{J_1} \vert  }
Q'_2 \bigg( {\mathcal B}^l_{\vert I  \vert} (x_I) . Q'_2 \Big({\mathcal B}^{m-l}_{\vert J_1 \vert} \big(Q_{\alpha,1}(x_{J_1})\big). {\mathcal B}^k_{\vert J_2 \vert} (x_{J_2})\Big)\bigg) \Bigg) \\
& = \sum^{\bullet}_{ I\sqcup J=[1;n] \atop \vert I\vert \geq 2} \varepsilon_x (I,J) \    (-1)^{\vert x_I \vert} (m+1) Q'_2 \Big({\mathcal B}_{\vert I \vert}^{m+1} \big( Q_{\alpha , 1} (x_I)\big) . {\mathcal B}_{\vert J \vert }^k (x_J)\Big)\\
& + \sum_{l=0}^m \sum^{\bullet}_{ I\sqcup J=[1;n] \atop \vert I \vert \geq 2 , \vert J\vert \geq 2} \varepsilon_x (I,J) (-1)^{\vert x_I  \vert } \      Q'_2 \bigg( {\mathcal B}^{m-l}_{\vert I \vert} \big(Q_{\alpha,1}(x_I)\big) . \sum^{\bullet}_{J_1 \sqcup J_2 = J    } \varepsilon_x (J_1,J_2) (-1)^{\vert x_{J_1} \vert } (\vert J_1 \vert -1)
 Q'_2 \Big( {\mathcal B}^l_{\vert J_1  \vert} (x_{J_1}) .  {\mathcal B}^k_{\vert J_2 \vert} (x_{J_2})\Big)\bigg)
  \end{flalign*}

\noindent et

\begin{flalign*}  
 & \sum_{l=0}^m \Bigg( \sum^{\bullet}_{ I\sqcup J=[1;n] \atop \vert I\vert \geq 3} \varepsilon_x (I,J) \    (-1)^{\vert x_I \vert +1 } (\vert I \vert -2) \sum^{\bullet}_{I_1 \sqcup I_2 = I \atop \vert I_1 \vert \geq 2} \varepsilon_x (I_1,I_2) 
Q'_2 \bigg( Q'_{\beta ,2} \Big({\mathcal B}^l_{\vert I_1  \vert} (x_{I_1}) \otimes {\mathcal B}^{m-l}_{\vert I_2 \vert} (x_{I_2})\Big). {\mathcal B}^k_{\vert J \vert} (x_{J})\bigg) \\ 
& + \sum^{\bullet}_{ I\sqcup J=[1;n] \atop \vert J\vert \geq 2} \varepsilon_x (I,J) \      (\vert I \vert -1) \sum^{\bullet}_{J_1 \sqcup J_2 = J \atop \vert J_1 \vert \geq 2} \varepsilon_x (J_1,J_2) 
Q'_2 \bigg( {\mathcal B}^l_{\vert I  \vert} (x_I) . Q'_{\beta ,2} \Big({\mathcal B}^{m-l}_{\vert J_1 \vert} (x_{J_1})\otimes {\mathcal B}^k_{\vert J_2 \vert} (x_{J_2})\Big)\bigg) \Bigg) \\
& = \sum^{\bullet}_{ I\sqcup J=[1;n] \atop \vert I\vert \geq 2} \varepsilon_x (I,J) \     (m+1) Q'_{\beta ,2} \Big({\mathcal B}_{\vert I \vert}^{m+1}  (x_I) \otimes {\mathcal B}_{\vert J \vert }^k (x_J)\Big)\\
& + \sum_{l=0}^m \sum^{\bullet}_{ I\sqcup J=[1;n] \atop \vert I \vert \geq 2 ,\vert J\vert \geq 2} \varepsilon_x (I,J) (-1)^{\vert x_I  \vert } \      Q'_2 \bigg( Q'_{\beta ,1} \big({\mathcal B}^{m-l}_{\vert I \vert} (x_I)\big) . \sum^{\bullet}_{J_1 \sqcup J_2 = J    } \varepsilon_x (J_1,J_2) (-1)^{\vert x_{J_1} \vert } (\vert J_1 \vert -1)
 Q'_2 \Big( {\mathcal B}^l_{\vert J_1  \vert} (x_{J_1}) .  {\mathcal B}^k_{\vert J_2 \vert} (x_{J_2})\Big)\bigg)
  \end{flalign*}

\end{lem}

\begin{preuve} (du lemme \ref{lem:sommes}).
 Pour expliciter les calculs,  on revient aux sommes respectueuses sur trois  parties de $ \{1,...,n\}$. Il faut tout simplement remarquer que si on casse le premier sous-ensemble 
 d'une somme respectueuse on obtient une somme respectueuse, mais que si on s\'epare le deuxi\` eme sous-ensemble, on obtient deux sommes dont l'une est respectueuse mais l'autre ne l'est pas.
 Mais cette derni\` ere le devient si on permute les deux premiers sous-ensembles, ce qui revient \`a multiplier par un signe.  
Le calcul se fait ainsi~:

\begin{flalign*}  
 &  \sum^{\bullet}_{ I\sqcup J=[1;n] \atop \vert I\vert \geq 3} \varepsilon_x (I,J) \    (-1)^{\vert x_I \vert +1 } (\vert I \vert -2) \sum^{\bullet}_{I_1 \sqcup I_2 = I \atop \vert I_1 \vert \geq 2} \varepsilon_x (I_1,I_2) (-1)^{\vert x_{I_1} \vert  }
Q'_2 \bigg( Q'_2 \Big({\mathcal B}^l_{\vert I_1  \vert} \big(Q_{\alpha,1}(x_{I_1}\big)\Big) . {\mathcal B}^{m-l}_{\vert I_2 \vert} (x_{I_2})). {\mathcal B}^k_{\vert J \vert} (x_{J})\bigg) \\ 
& + \sum^{\bullet}_{ I\sqcup J=[1;n] \atop \vert J\vert \geq 2} \varepsilon_x (I,J) \      (\vert I \vert -1) \sum^{\bullet}_{J_1 \sqcup J_2 = J \atop \vert J_1 \vert \geq 2} \varepsilon_x (J_1,J_2) (-1)^{\vert x_{J_1} \vert  }
Q'_2 \bigg( {\mathcal B}^{m-l}_{\vert I  \vert} (x_I) . Q'_2 \Big({\mathcal B}^l_{\vert J_1 \vert} \big(Q_{\alpha,1}(x_{J_1})\big). {\mathcal B}^k_{\vert J_2 \vert} (x_{J_2})\Big)\bigg) \\ 
 & = \sum^{\bullet}_{ I\sqcup J=[1;n] \atop \vert I\vert \geq 3} \varepsilon_x (I,J) \    (-1)^{\vert x_I \vert +1 } (\vert I \vert -2) \sum^{\bullet}_{I_1 \sqcup I_2 = I \atop \vert I_1 \vert \geq 2} \varepsilon_x (I_1,I_2) (-1)^{\vert x_{I_1} \vert  }
Q'_2 \bigg( Q'_2 \Big({\mathcal B}^l_{\vert I_1  \vert} \big(Q_{\alpha,1}(x_{I_1}\big)\Big) . {\mathcal B}^{m-l}_{\vert I_2 \vert} (x_{I_2})). {\mathcal B}^k_{\vert J \vert} (x_{J})\bigg) \\ 
& + \sum^{\bullet}_{ I\sqcup J=[1;n] \atop \vert J\vert \geq 2} \varepsilon_x (I,J) \     (\vert I \vert -1) \sum^{\bullet}_{J_1 \sqcup J_2 = J \atop \vert J_1 \vert \geq 2} \varepsilon_x (J_1,J_2) (-1)^{\vert x_{J_1} \vert +\vert x_I  \vert}
Q'_2 \bigg( Q'_2 \Big( {\mathcal B}^{m-l}_{\vert I  \vert} (x_I) . {\mathcal B}^l_{\vert J_1 \vert} \big(Q_{\alpha,1}(x_{J_1})\big)\Big). {\mathcal B}^k_{\vert J_2 \vert} (x_{J_2})\bigg) \\
& + \sum^{\bullet}_{ I\sqcup J=[1;n] \atop \vert J\vert \geq 2} \varepsilon_x (I,J) \      (\vert I \vert -1) \sum^{\bullet}_{J_1 \sqcup J_2 = J \atop \vert J_1 \vert \geq 2} \varepsilon_x (J_1,J_2) (-1)^{\vert x_{J_1} \vert + \vert x_I \vert ( \vert x_{J_1} \vert +1)  }
Q'_2 \bigg( {\mathcal B}^l_{\vert J_1 \vert} \big(Q_{\alpha,1}(x_{J_1})\big) . \Big ({\mathcal B}^{m-l}_{\vert I  \vert} (x_I) .  {\mathcal B}^k_{\vert J_2 \vert} (x_{J_2})\Big)\bigg)
\end{flalign*}  

\begin{flalign*}  
& = \sum^{\bullet}_{ K_1\sqcup K_2 \sqcup K_3 =[1;n] \atop \vert K_1\vert + \vert K_2 \vert \geq 3} \varepsilon_x (K_1,K_2,K_3) \    (-1)^{\vert x_{K_2} \vert +1 } (\vert K_1 \vert + \vert K_2 \vert -2) 
Q'_2 \bigg( Q'_2 \Big({\mathcal B}^l_{\vert K_1  \vert} \big(Q_{\alpha,1}(x_{K_1}\big)\Big) . {\mathcal B}^{m-l}_{\vert K_2 \vert} (x_{K_2})). {\mathcal B}^k_{\vert K_3 \vert} (x_{K_3})\bigg) \\ 
& + \sum^{\bullet}_{ K_1\sqcup K_2 \sqcup K_3=[1;n] } \varepsilon_x (K_1,K_2,K_3) \     (\vert K_1 \vert -1)  (-1)^{\vert x_{K_1} \vert +\vert x_{K_2}  \vert}
Q'_2 \bigg( Q'_2 \Big( {\mathcal B}^{m-l}_{\vert K_1  \vert} (x_{K_1}) . {\mathcal B}^l_{\vert K_2 \vert} \big(Q_{\alpha,1}(x_{K_2})\big)\Big). {\mathcal B}^k_{\vert K_3 \vert} (x_{K_3})\bigg) \\
& + \sum^{\bullet}_{ K_1\sqcup K_2 \sqcup K_3=[1;n] } \varepsilon_x (K_1,K_2,K_3) (-1)^{ \vert x_{K_1} \vert \vert x_{K_2} \vert }\     (\vert K_2 \vert -1)  (-1)^{\vert x_{K_1} \vert +\vert x_{K_2}  \vert} \\
& \hspace{3cm} Q'_2 \bigg( Q'_2 \Big( {\mathcal B}^{m-l}_{\vert K_2  \vert} (x_{K_2}) . {\mathcal B}^l_{\vert K_1 \vert} \big(Q_{\alpha,1}(x_{K_1})\big)\Big). {\mathcal B}^k_{\vert K_3 \vert} (x_{K_3})\bigg) \\
& + \sum^{\bullet}_{ K_1\sqcup K_2 \sqcup K_3=[1;n]} \varepsilon_x (K_1,K_2,K_3) \      (\vert K_1 \vert -1)  (-1)^{\vert x_{K_1} \vert + \vert x_{K_2} \vert + \ \vert x_{K_1} \vert \vert x_{K_2} \vert }
Q'_2 \bigg( {\mathcal B}^l_{\vert K_2 \vert} \big(Q_{\alpha,1}(x_{K_2})\big) . \Big ({\mathcal B}^{m-l}_{\vert K_1  \vert} (x_{K_1}) .  {\mathcal B}^k_{\vert K_3 \vert} (x_{K_3})\Big)\bigg)  \\
& + \sum^{\bullet}_{ K_1\sqcup K_2 \sqcup K_3=[1;n] } \varepsilon_x (K_1,K_2,K_3)  x_{K_1} \vert  \      (\vert K_2 \vert -1)  (-1)^{\vert x_{K_1} \vert + \vert x_{K_2} \vert  }
Q'_2 \bigg( {\mathcal B}^l_{\vert K_1 \vert} \big(Q_{\alpha,1}(x_{K_1})\big) . \Big ({\mathcal B}^{m-l}_{\vert K_2  \vert} (x_{K_2}) .  {\mathcal B}^k_{\vert K_3 \vert} (x_{K_3})\Big)\bigg) \\
 & = \sum^{\bullet}_{ I\sqcup J=[1;n] \atop \vert I\vert \geq 2} \varepsilon_x (I,J) \    (-1)^{\vert x_I \vert   }  Q'_2 \bigg[ \sum^{\bullet}_{I_1 \sqcup I_2 = I  } \varepsilon_x (I_1,I_2)
 (\vert I_1 \vert -1) \\
 & \hspace{3cm} \bigg( (-1)^{\vert x_{I_1} \vert  +1}
   Q'_2 \Big({\mathcal B}^l_{\vert I_1  \vert} \big(Q_{\alpha,1}(x_{I_1}\big)  . {\mathcal B}^{m-l}_{\vert I_2 \vert} (x_{I_2})\Big)+Q'_2 \Big({\mathcal B}^{m-l}_{\vert I_1  \vert} (x_{I_1} 
   . {\mathcal B}^l_{\vert I_2 \vert} \big(Q_{\alpha,1}(x_{I_2})\big)\Big)\bigg). {\mathcal B}^k_{\vert J \vert} (x_{J})\bigg] \\ 
& + \sum^{\bullet}_{ I\sqcup J=[1;n] \atop \vert J\vert \geq 2} \varepsilon_x (I,J) (-1)^{\vert x_I  \vert } \      Q'_2 \bigg( {\mathcal B}^l_{\vert I \vert} \big(Q_{\alpha,1}(x_I)\big). \sum^{\bullet}_{J_1 \sqcup J_2 = J    } \varepsilon_x (J_1,J_2) (-1)^{\vert x_{J_1} \vert } (\vert J_1 \vert -1)
 Q'_2 \Big( {\mathcal B}^{m-l}_{\vert J_1  \vert} (x_{J_1}) .  {\mathcal B}^k_{\vert J_2 \vert} (x_{J_2})\Big)\bigg)   
 \end{flalign*}

\noindent En faisant la somme sur $l$ et en utilisant le lemme \ref{lem:sommes}, on d\'eduit le r\'esultat. La preuve de la deuxi\` eme formule est analogue puisque~:

$$ Q'_{\beta ,2} \Big({\mathcal B}^l_{\vert I_1  \vert} (x_{I_1}) \otimes {\mathcal B}^{m-l}_{\vert I_2 \vert} (x_{I_2})\Big) = 
(-1)^{\vert x_{I_1} \vert } Q'_2 \Big(Q'_{\beta ,1}\big({\mathcal B}^l_{\vert I_1  \vert} (x_{I_1}) \big). {\mathcal B}^{m-l}_{\vert I_2 \vert} (x_{I_2})\Big) $$

\noindent et~:

\begin{flalign*} & (-1)^{\vert x_{I_1} \vert } Q'_{\beta ,1} \Big( Q'_2 \big({\mathcal B}^l_{\vert I_1  \vert} (x_{I_1}) . 
{\mathcal B}^{m-l}_{\vert I_2 \vert} (x_{I_2})\big)\Big)  \\ & = (-1)^{\vert x_{I_1} \vert +1 }  
Q'_2 \Big(Q'_{\beta ,1}\big({\mathcal B}^l_{\vert I_1  \vert} (x_{I_1}) \big). {\mathcal B}^{m-l}_{\vert I_2 \vert} (x_{I_2})\Big) +
Q'_2 \Big({\mathcal B}^l_{\vert I_1  \vert} (x_{I_1}) . Q'_{\beta ,1}\big({\mathcal B}^{m-l}_{\vert I_2 \vert} (x_{I_2})\big)\Big).
\end{flalign*}

\end{preuve}

Montrons maintenant la proposition \ref{prop:Cnj}.
\begin{preuve} On va d\'emontrer la r\'esultat par r\'ecurrence. La remarque \ref{rmque:Cn0} donne l'initialisation pour $j=0$.
Supposons que la propri\'et\'e \'enonc\'ee dans la proposition \ref{prop:Cnj} est vraie jusqu'\`a l'ordre $j-1$. On a alors~:

\begin{flalign*} {\mathcal C}_n^j (x_1\otimes ... \otimes x_n)& =  
  \frac{1}{j} \sum^{\bullet}_{ I\sqcup J=[1;n] } \varepsilon_x (I,J) \  \  (-1)^{\vert x_I \vert +1 } (\vert I \vert -1) \sum_{k=0}^{j-1} 
Q'_2 ({\mathcal B}^k_{\vert I  \vert} (Q_{\alpha,1}(x_I)) . {\mathcal B}^{j-k-1}_{\vert J \vert} (x_J)) \\ 
& +
\frac{1}{j} \sum^{\bullet}_{ I\sqcup J=[1;n] } \varepsilon_x (I,J) \  \   (\vert I \vert -1) \sum_{k=0}^{j-1} 
Q'_2 ({\mathcal B}^k_{\vert I  \vert} (x_I) .{\mathcal B}^{j-k-1}_{\vert J \vert } (Q_{\alpha,1}(x_J)))\\
& + 
  \frac{1}{j} \sum^{\bullet}_{ I\sqcup J=[1;n] \atop \vert I \vert \geq 2 } \varepsilon_x (I,J) \  \  (-1)^{\vert x_I \vert +1 } (\vert I \vert -2) \sum_{k=0}^{j-1} 
Q'_2 ({\mathcal B}^k_{\vert I  \vert-1} (Q_{\alpha,2}(x_I)) . {\mathcal B}^{j-k-1}_{\vert J \vert } (x_J)) \\ 
& +
\frac{1}{j} \sum^{\bullet}_{ I\sqcup J=[1;n] \atop \vert J \vert \geq 2} \varepsilon_x (I,J) \  \   (\vert I \vert -1) \sum_{k=0}^{j-1} 
Q'_2 ({\mathcal B}^k_{\vert I  \vert} (x_I) .{\mathcal B}^{j-k-1}_{\vert J \vert -1} (Q_{\alpha,2}(x_J)))\\
& +
\frac{1}{j} \sum^{\bullet}_{ I\sqcup J=[1;n] } \varepsilon_x (I,J) \  (-1)^{\vert x_I \vert  } \   (\vert I \vert -1) \sum_{k=0}^{j-1} 
Q'_2 (Q'_{\beta ,1}({\mathcal B}^k_{\vert I  \vert} (x_I)) .{\mathcal B}^{j-k-1}_{\vert J \vert } (x_J)) \\
& -
\frac{1}{j} \sum^{\bullet}_{ I\sqcup J=[1;n] } \varepsilon_x (I,J) \  \   (\vert I \vert -1) \sum_{k=0}^{j-1} 
Q'_2 ({\mathcal B}^k_{\vert I  \vert} (x_I) . Q'_{\beta ,1}({\mathcal B}^{j-k-1}_{\vert J \vert } (x_J))) \\
&= 
\frac{1}{j} \sum^{\bullet}_{ I\sqcup J=[1;n] \atop \vert I \vert \geq 2} \varepsilon_x (I,J) \       \sum_{k=0}^{j-1} 
Q'_{\beta ,2} ({\mathcal B}^k_{\vert I  \vert} (x_I) \otimes {\mathcal B}^{j-k-1}_{\vert J \vert } (x_J))
\\
&+
\frac{1}{j} \sum^{\bullet}_{ I\sqcup J=[1;n] \atop \vert I \vert \geq 2} \varepsilon_x (I,J) \  \  (-1)^{\vert x_I \vert +1 } \ \sum_{k=0}^{j-1} 
Q'_2 ({\mathcal B}^k_{\vert I  \vert} (Q_{\alpha,1}(x_I)) . {\mathcal B}^{j-k-1}_{\vert J \vert} (x_J)) \\ 
&+
\frac{1}{j} \sum^{\bullet}_{ I\sqcup J=[1;n] \atop \vert I \vert \geq 3} \varepsilon_x (I,J) \  \  (-1)^{\vert x_I \vert +1 } \ (\vert I \vert -2) \  \sum_{k=0}^{j-1} 
Q'_2 ({\mathcal C}^{j-k-1}_{\vert I  \vert} (x_I) . {\mathcal B}^k_{\vert J \vert} (x_J)) \\ 
&+
\frac{1}{j} \sum^{\bullet}_{ I\sqcup J=[1;n] \atop \vert J \vert \geq 2} \varepsilon_x (I,J) \  \   (\vert I \vert -1) \ \sum_{k=0}^{j-1} 
Q'_2 ({\mathcal B}^{j-k-1}_{\vert I  \vert} (x_I) . {\mathcal C}^k_{\vert J \vert} (x_J)) \\
& =
\frac{1}{j} \sum^{\bullet}_{ I\sqcup J=[1;n] \atop \vert I \vert \geq 2} \varepsilon_x (I,J) \       \sum_{k=0}^{j-1} 
Q'_{\beta ,2} ({\mathcal B}^k_{\vert I  \vert} (x_I) \otimes {\mathcal B}^{j-k-1}_{\vert J \vert } (x_J))
 + \frac{1}{j} R_n^j + V_n^j
\end{flalign*}
o\`u~:
\begin{flalign*}
   R_n^j   & = \frac{1}{j} \sum^{\bullet}_{ I\sqcup J=[1;n] \atop \vert I \vert \geq 3} \varepsilon_x (I,J) \  \  (-1)^{\vert x_I \vert +1 } \ (\vert I \vert -2) \  \sum_{k=0}^{j-1} 
Q'_2 ({\mathcal C}^{j-k-1}_{\vert I  \vert} (x_I) . {\mathcal B}^k_{\vert J \vert} (x_J)) \\
& +
\frac{1}{j} \sum^{\bullet}_{ I\sqcup J=[1;n] \atop \vert J \vert \geq 2} \varepsilon_x (I,J) \  \   (\vert I \vert -1) \ \sum_{k=0}^{j-1} 
Q'_2 ({\mathcal B}^{j-k-1}_{\vert I  \vert} (x_I) . {\mathcal C}^k_{\vert J \vert} (x_J))
 \end{flalign*}

\noindent Il s'agit de voir qu'une fois \'ecrits seulement \` a l'aide des $ {\mathcal B}_s^0$, les termes $ {\mathcal C}_q^m $ ou $V_q^m$ ( $ m \geq 1$) sont la somme de deux quantit\'es "homog\` enes"
vis-\`a-vis du nombre d'occurrence de $ Q'_2$. 
Plus exactement, la premi\`ere contient $ m $ occurrences de $Q'_2$ et l'autre quantit\'e en contient $ m+1$. Ceci n'est pas vrai pour ${ \cal C}_q^0$, la partie ne contenant aucune occurrence de $ Q'_2$ est nulle.

\noindent La partie de $ V_n^j $ contenant $j$ occurrences de $ Q'_2$ est donn\'ee par~:
\begin{flalign*}
 &\frac{1}{j} \sum^{\bullet}_{ I\sqcup J=[1;n] \atop \vert I \vert \geq 3} \varepsilon_x (I,J) \  \  (-1)^{\vert x_I \vert +1 } \ (\vert I \vert -2) \  \sum_{k=0}^{j-2} 
Q'_2  \Big(\biggl[\sum^{\bullet}_{I_1 \sqcup I_2 = I \atop \vert I_1 \vert \geq 2} \varepsilon_x (I_1,I_2)  
\\
&\sum_{l=0}^{j-k-2 }  
   \Bigl[Q'_{\beta ,2}\big( {\mathcal B}^{l}_{\vert I_1  \vert} (x_{I_1}) \otimes {\mathcal B}^{j-k-l-2}_{\vert I_2 \vert} (x_{I_2}) \big)-(-1)^{\vert x_{I_1} \vert} Q'_2 \Big({\mathcal B }^l_{\vert I_1 \vert }\big(Q_{\alpha , 1} (x_{I_1})\big) . {\mathcal B}^{j-k-l-2}_{\vert I_2 \vert} (x_{I_2})\Big)\Bigr]\biggr]. {\mathcal B}^k_{\vert J \vert} (x_J)\bigg) \\ 
&+
\frac{1}{j} \sum^{\bullet}_{ I\sqcup J=[1;n] \atop \vert J \vert \geq 2} \varepsilon_x (I,J) \    (\vert I \vert -1) \ \sum_{l=1}^{j-1} 
Q'_2 \bigg({\mathcal B}^{j-l-1}_{\vert I  \vert} (x_I) . \biggl[ \sum^{\bullet}_{J_1 \sqcup J_2 = I \atop \vert J_1 \vert \geq 2} \varepsilon_x (J_1,J_2)  \\
& \sum_{k=0}^{l-1 } \Bigl[Q'_{\beta , 2} \big(  {\mathcal B}^{l-k-1}_{\vert J_1 \vert} (x_{J_1})\otimes {\mathcal B}^k_{\vert J_2 \vert } (x_{J_2})\big) -  (-1)^{\vert x_{J_1} \vert} Q'_2 \Big({\mathcal B }^{l-k-1}_{\vert J_1 \vert }\big(Q_{\alpha , 1} (x_{J_1})\big) . {\mathcal B}^k_{\vert J_2 \vert} (x_{J_2})\Big)\Bigr]\biggr] \bigg) \\
&= \frac{1}{j} \sum_{k=0}^{j-2} \sum_{l=0}^{j-k-2}\biggl[
\sum^{\bullet}_{ I\sqcup J=[1;n] \atop \vert I \vert \geq 3} \varepsilon_x (I,J) \     (-1)^{\vert x_I \vert +1 } \ (\vert I \vert -2)  \   \sum^{\bullet}_{I_1 \sqcup I_2 = I \atop \vert I_1 \vert \geq 2} \varepsilon_x (I_1,I_2)   \\
& Q'_2  \Big(\Bigl[ Q'_{\beta ,2}\big( {\mathcal B}^{l}_{\vert I_1  \vert} (x_{I_1}) \otimes {\mathcal B}^{j-k-l-2}_{\vert I_2 \vert} (x_{I_2}) \big)-(-1)^{\vert x_{I_1} \vert} Q'_2 \Big({\mathcal B }^l_{\vert I_1 \vert }\big(Q_{\alpha , 1} (x_{I_1})\big) . {\mathcal B}^{j-k-l-2}_{\vert I_2 \vert} (x_{I_2})\Big)\Bigr]. {\mathcal B}^k_{\vert J \vert} (x_J)\Big) \\ 
&+
 \sum^{\bullet}_{ I\sqcup J=[1;n] \atop \vert J \vert \geq 2} \varepsilon_x (I,J) \    (\vert I \vert -1)  \sum^{\bullet}_{J_1 \sqcup J_2 = I \atop \vert J_1 \vert \geq 2} \varepsilon_x (J_1,J_2)  
Q'_2 \bigg({\mathcal B}^{j-k-l-2}_{\vert I  \vert} (x_I) .     \\
&    \Bigl[Q'_{\beta , 2} \big(  {\mathcal B}^{l}_{\vert J_1 \vert} (x_{J_1})\otimes {\mathcal B}^k_{\vert J_2 \vert } (x_{J_2})\big) -  (-1)^{\vert x_{J_1} \vert} Q'_2 \Big({\mathcal B }^{l}_{\vert J_1 \vert }\big(Q_{\alpha , 1} (x_{J_1})\big) . {\mathcal B}^k_{\vert J_2 \vert} (x_{J_2})\Big)\Bigr] \bigg) \biggr]
\\ 
& =  \frac{1}{j} \sum_{k=0}^{j-2} \sum^{\bullet}_{ I\sqcup J=[1;n] \atop \vert I \vert \geq 2} \varepsilon_x (I,J) (j-k-1) 
\Bigl[ Q'_{\beta ,2} \big( {\mathcal B}^{j-k-1}_{\vert I  \vert} (x_I) \otimes {\mathcal B}^k_{\vert J  \vert} (x_J)\big) - (-1)^{\vert x_I \vert }  
  Q'_2 \Big( {\mathcal B}^{j-k-1}_{\vert I  \vert} \big( Q_{\alpha ,1 } (x_I) \big). 
  {\mathcal B}^k_{\vert J  \vert} (x_J)\Big)\Bigr]  \\
  + & \frac{1}{j} \sum_{k=0}^{j-2} \sum_{l=0}^{j-k-2} 
\sum^{\bullet}_{ I\sqcup J=[1;n] \atop \vert J \vert \geq 2} \varepsilon_x (I,J) (-1)^{\vert x_I  \vert } Q'_2 \Big( \Bigr[Q'_{\beta ,1} \big( {\mathcal B}^l_{\vert I \vert } (x_I)\big) -  {\mathcal B}^l_{\vert I \vert } \big(Q_{\alpha ,1} (x_I)\big)\Bigr].\\
 & \sum^{\bullet}_{J_1 \sqcup J_2 = J } \varepsilon_x (J_1,J_2) 
   (-1)^{ + \vert x_{J_1}  }  (\vert J_1 \vert -1)  
       Q'_2 \big( {\mathcal B}^{j-k-l-2}_{\vert I \vert } (x_{J_1}).{\mathcal B}^k_{\vert J_2 \vert } (x_{J_2})\big)\Big)  
\end{flalign*}
\begin{flalign*}   
 & =  \frac{1}{j} \sum_{k=0}^{j-2} \sum^{\bullet}_{ I\sqcup J=[1;n] \atop \vert I \vert \geq 2} \varepsilon_x (I,J) (j-k-1)\Bigl[ Q'_{\beta ,2} \big( {\mathcal B}^{j-k-1}_{\vert I  \vert} (x_I) \otimes {\mathcal B}^k_{\vert J  \vert} (x_J)\big) - Q'_2 \big( {\mathcal B}^{j-k-1}_{\vert I  \vert} \big( Q_{\alpha , 1}(x_I)\big) . {\mathcal B}^{k}_{\vert J  \vert} (x_J)\big)\Bigr] \\ 
   &+ \frac{1}{j} \sum_{l=0}^{j-2} \sum^{\bullet}_{ I\sqcup J=[1;n] \atop \vert I \vert \geq 2} \varepsilon_x (I,J) (j-l-1)\Bigl[ Q'_{\beta ,2} \big( {\mathcal B}^{l}_{\vert I  \vert} (x_I) \otimes {\mathcal B}^{j-l-1}_{\vert J  \vert} (x_J)\big) - Q'_2 \big( {\mathcal B}^{l}_{\vert I  \vert} \big( Q_{\alpha , 1}(x_I)\big) . {\mathcal B}^{j-l-1}_{\vert J  \vert} (x_J)\big)\Bigr]\\
   & = \frac{1}{j} \sum_{k=1}^{j-1} \sum^{\bullet}_{ I\sqcup J=[1;n] \atop \vert I \vert \geq 2} \varepsilon_x (I,J) k \Bigl[ Q'_{\beta ,2} \big( {\mathcal B}^{k}_{\vert I  \vert} (x_I) \otimes {\mathcal B}^{j-k-1}_{\vert J  \vert} (x_J)\big) - Q'_2 \big( {\mathcal B}^{k}_{\vert I  \vert} \big( Q_{\alpha , 1}(x_I)\big) . {\mathcal B}^{j-k-1}_{\vert J  \vert} (x_J)\big)\Bigr]\\   
   &+ \frac{1}{j} \sum_{k=0}^{j-2} \sum^{\bullet}_{ I\sqcup J=[1;n] \atop \vert I \vert \geq 2} \varepsilon_x (I,J) (j-k-1) \Bigl[ Q'_{\beta ,2} \big( {\mathcal B}^{k}_{\vert I  \vert} (x_I) \otimes {\mathcal B}^{j-k-1}_{\vert J  \vert} (x_J)\big)- Q'_2 \big( {\mathcal B}^{k}_{\vert I  \vert} \big( Q_{\alpha , 1}(x_I)\big) . {\mathcal B}^{j-k-1}_{\vert J  \vert} (x_J)\big)\Bigr]\\
  & = \frac{j-1}{j}  \sum^{\bullet}_{ I\sqcup J=[1;n] \atop \vert I \vert \geq 2} \varepsilon_x (I,J)  \sum_{k=0}^{j-1} Q'_{\beta ,2} \big( {\mathcal B}^k_{\vert I  \vert} (x_I) \otimes {\mathcal B}^{j-k-1}_{\vert J  \vert} (x_J)\big) + \frac{j-1}{j} R_n^j
\end{flalign*}

Pour les termes contenant $ (j+1)$ occurrences de $Q'_2$, on commence par remarquer que 

$$\sum^{\bullet}_{ I\sqcup J=[1;n] \atop \vert I \vert = 1 }  \varepsilon_x (I,J) Q'_{\beta,2} \big({\mathcal B}_1^k (x_I) \otimes {\mathcal B}_{n-1}^{m-k} (x_J) \big) = \sum^{\bullet}_{ I\sqcup J=[1;n] \atop \vert I \vert = 1 }  \varepsilon_x (I,J) (-1)^{\vert x_I \vert } Q'_2 \big({\mathcal B}_1^k \big( Q_{\alpha , 1}(x_I)\big)  {\mathcal B}_{n-1}^{m-k} (x_J) \big). $$

\noindent On peut donc  rajouter ces termes  \` a $ - R_n^{m+1} $ et on supprime la condition      $ \vert I \vert \geq 2$. Il n'est pas difficile de voir que dans les \'egalit\'es  du lemme 3, la suppression des conditions 
$ \vert I_1 \vert \geq 2 $ et $ \vert J_1 \vert \geq 2 $ dans le premier membre se traduit par la suppression de la condition 
$ \vert  I  \vert \geq 2 $ dans le second membre. Les termes contenant $ (j+1)$  occurrences de $ Q'_2$ seront alors~:
\begin{flalign*}
& \frac{1}{j} \sum_{k=0}^{j-1} \sum_{l=0}^{j-k-1}\biggl[
\sum^{\bullet}_{ I\sqcup J=[1;n] \atop \vert I \vert \geq 3} \varepsilon_x (I,J) \     (-1)^{\vert x_I \vert +1 } \ (\vert I \vert -2)  \   \sum^{\bullet}_{I_1 \sqcup I_2 = I \atop \vert I_1 \vert \geq 2} \varepsilon_x (I_1,I_2)   \\
& Q'_2  \Big((-1)^{\vert x_{I_1} \vert} Q'_2 \Big({\mathcal B }^l_{\vert I_1 \vert }\big(Q_{\alpha , 1} (x_{I_1})\big) . {\mathcal B}^{j-k-l-1}_{\vert I_2 \vert} (x_{I_2})\Big). {\mathcal B}^k_{\vert J \vert} (x_J)\Big) \\ 
&+
 \sum^{\bullet}_{ I\sqcup J=[1;n] \atop \vert J \vert \geq 2} \varepsilon_x (I,J) \    (\vert I \vert -1)  \sum^{\bullet}_{J_1 \sqcup J_2 = I \atop \vert J_1 \vert \geq 2} \varepsilon_x (J_1,J_2)  
Q'_2 \bigg({\mathcal B}^{j-k-l-1}_{\vert I  \vert} (x_I) .     \\
&      (-1)^{\vert x_{J_1} \vert} Q'_2 \Big({\mathcal B }^{l}_{\vert J_1 \vert }\big(Q_{\alpha , 1} (x_{J_1})\big) . {\mathcal B}^k_{\vert J_2 \vert} (x_{J_2})\Big)  \bigg) \biggr]\\
\end{flalign*}

\noindent Par les m\^ emes transformations que ci-dessus, ces termes sont \'egaux \` a
\begin{flalign*}
&  \frac{1}{j} \sum_{k=0}^{j-1} \sum^{\bullet}_{ I\sqcup J=[1;n] } \varepsilon_x (I,J) (j-k ) Q'_2 \big( {\mathcal B}^{k}_{\vert I  \vert} \big( Q_{\alpha , 1}(x_I)\big) . {\mathcal B}^{j-k }_{\vert J  \vert} (x_J)\big)\Bigr]\\
  & +\frac{1}{j} \sum_{k=1}^{j } \sum^{\bullet}_{ I\sqcup J=[1;n] } \varepsilon_x (I,J) ( k ) Q'_2 \big( {\mathcal B}^{k}_{\vert I  \vert} \big( Q_{\alpha , 1}(x_I)\big) . {\mathcal B}^{j-k }_{\vert J  \vert} (x_J)\big)\Bigr] \\
  & =\sum^{\bullet}_{ I\sqcup J=[1;n] \atop \vert I \vert = 1 }  \varepsilon_x (I,J) \sum_{k=0}^j Q'_{\beta,2} \big({\mathcal B}_1^k (x_I) \otimes {\mathcal B}_{n-1}^{j-k} (x_J) \big) - R_n^{j+1} (x_1\otimes ... \otimes x_n).
  \end{flalign*}
  \end{preuve}
  
D\'emontrons maintenant le th\'eor\`eme \ref{theo:principal}.
\begin{preuve}
Donnons-nous des applications ${\mathcal B}_n : T^{n }{\mathfrak g} \to  {\mathfrak g}'$ de degr\'e $0$ qui v\'erifient les relations suivantes:
\begin{eqnarray*}   {\mathcal B}_n \circ Q_{\alpha,1}  \ \  (x_1 \otimes \cdots \otimes x_n )  + {\mathcal B}_{n-1}   \circ Q_{\alpha,2}  \ \  (x_1 \otimes \cdots \otimes x_n )  &=& \\
 Q_{\beta,1}' \circ  {\mathcal B}_n   \ \ 
(x_1 \otimes \cdots \otimes x_n )  &  & \\
+ \sum^{\bullet}_{ I\sqcup J=[1;n] }  \varepsilon_x (I,J) \ \ Q'_{\beta,2}
 \left(  {\mathcal B}_i (x_I)  \otimes  {\mathcal B}_j (x_J)  \right) & &
\end{eqnarray*}
Par la proposition \ref{prop:toutsetend}, les applications $({\mathcal B}_n)_{ n\geq 1} $ s'\'etendent en un morphisme bien fait. 

\noindent
Par d\'efinition des ${\mathcal C}_n^j$, pour tout entier $ n \geq 2$ et tout $x_1, \dots,x_n \in {\mathfrak g}$~:
\begin{equation} 
 \label{eq:sommeBij}
 \sum_{j \geq 0} {\mathcal C}_n^j  = \sum_{j ^\geq 0}  \left( B_n^j \circ Q_{\alpha,1}  + B_{n-1}^j  \circ Q_{\alpha,1} - Q_{\beta,1}' \circ  B_n^j 
\right)(x_1 \otimes \cdots \otimes x_n )  .
\end{equation}
Par les propositions \ref{prop:Cn0} et \ref{prop:Cnj}, on a~:
\begin{eqnarray*} 
 \sum_{j \geq 0} {\mathcal C}_n^j &=&  {\mathcal C}_ n^0 + \sum_{j \geq 1} {\mathcal C}_n^j  \\
 &=&  \sum^{\bullet}_{ I\sqcup J=[1;n] \atop \vert I \vert = 1 }  \varepsilon_x (I,J) Q'_{\beta,2} \big({\mathcal B}_1^0 (x_I) \otimes {\mathcal B}_{n-1}^0 (x_J) \big) -R_n^1 \\
 & & +   \sum_{j \geq 1} \sum^{\bullet}_{ I\sqcup J=[1;n] \atop \vert I \vert \geq 2 }
 \varepsilon_x (I,J) \sum_{k=0}^{j-1} Q'_{\beta,2} \big({\mathcal B}_{\vert I \vert }^k (x_I) \otimes {\mathcal B}_{\vert J \vert }^{j-k-1} (x_J) \big) \\ 
 & & +  \sum_{j \geq 1} \sum^{\bullet}_{ I\sqcup J=[1;n] \atop \vert I \vert = 1 }  \varepsilon_x (I,J)
  \sum_{j \geq 1}\sum_{k=0}^j  Q'_{\beta,2} \big({\mathcal B}_1^k (x_I) \otimes {\mathcal B}_{n-1}^{j-k} (x_J) \big)  \\
 & & +  \sum_{j \geq 1} \left( R_n^j (x_1 \otimes... \otimes x_n) -   R_n^{j+1} (x_1 \otimes... \otimes x_n) \right)
\end{eqnarray*}
Les termes $(R_n^j)_{j \geq 1}$  se simplifient deux \`a deux~:
\begin{eqnarray*} 
 \sum_{j \geq 0} {\mathcal C}_n^j &=& \sum^{\bullet}_{ I\sqcup J=[1;n] \atop \vert I \vert = 1 }  \varepsilon_x (I,J) Q'_{\beta,2}
 \big({\mathcal B}_1^0 (x_I) \otimes {\mathcal B}_{n-1}^0 (x_J) \big) \\
 & & +   \sum_{j \geq 1} \sum^{\bullet}_{ I\sqcup J=[1;n] \atop \vert I \vert \geq 2 }
 \varepsilon_x (I,J) \sum_{k=0}^{j-1} Q'_{\beta,2} \big({\mathcal B}_{\vert I \vert }^k (x_I) \otimes {\mathcal B}_{\vert J \vert }^{j-k-1} (x_J) \big) \\ 
 & & +  \sum_{j \geq 1} \sum^{\bullet}_{ I\sqcup J=[1;n] \atop \vert I \vert = 1 }  \varepsilon_x (I,J)
  \sum_{k=0}^j  Q'_{\beta,2} \big({\mathcal B}_1^k (x_I) \otimes {\mathcal B}_{n-1}^{j-k} (x_J) \big)  \\
\end{eqnarray*}
 Les sommes du c\^ot\'e droit se r\'eunissent en une seule somme par un simple jeu sur les indices~:
\begin{eqnarray*} 
 \sum_{j \geq 0} {\mathcal C}_n^j &=& \sum^{\bullet}_{ I\sqcup J=[1;n] }  \varepsilon_x (I,J) \ \ Q'_{\beta,2}
 \left( \left( \sum_{k \geq 0} {\mathcal B}_{|i|}^k (x_I) \right) \otimes  \left( \sum_{l \geq 0} {\mathcal B}_{|J|}^l (x_J) \right)  \right) \\ 
\end{eqnarray*} 
En comparant avec l'\'equation (\ref{eq:sommeBij}), on obtient pour tous $x_1, \dots,x_n \in {\mathfrak g}$~:
\begin{eqnarray*}   \left( \sum_{j ^\geq 0} B_n^j  \right) \circ Q_{\alpha,1}  \ \  (x_1 \otimes \cdots \otimes x_n )  + \left( \sum_{j ^\geq 0} B_{n-1}^j  \right)    \circ Q_{\alpha,2}  \ \  (x_1 \otimes \cdots \otimes x_n )  &=& \\
 Q_{\beta,1}' \circ \left( \sum_{j ^\geq 0} B_n^j  \right)   \ \ 
(x_1 \otimes \cdots \otimes x_n )  &  & \\
+ \sum^{\bullet}_{ I\sqcup J=[1;n] }  \varepsilon_x (I,J) \ \ Q'_{\beta,2}
 \left( \left( \sum_{k \geq 0} {\mathcal B}_{|i|}^k (x_I) \right) \otimes  \left( \sum_{l \geq 0} {\mathcal B}_{|J|}^l (x_J) \right)  \right) & &
\end{eqnarray*}
De la proposition \ref{prop:toutsetend}, il suit que la famille $ (\sum_{j ^\geq 0} B_n^j )_{ n \geq 1} $ est la famille des coefficients de Taylor
d'un $ {\rm Leib}_\infty$-morphisme bien fait.
\end{preuve}

\section{Applications et questions}  
\newcommand{\planck}{h}

Nous allons utiliser notre construction pour retrouver la formule due \`a Dominique Manchon \cite{Manchon}, qui relie -selon ses mots- le crochet de Poisson, le commutateur 
du star-produit, l'application tangente du morphisme de formalit\'e de Kontsevich, et sa d\'eriv\'ee seconde.

On prend les notations de \cite{AMM}-\cite{Kont}-\cite{Manchon}, qui sont d\'esormais classiques. Soit $M$ une vari\'et\'e diff\'erentielle, et
soit $\gamma \in T_{poly}(M)[[ \hbar]] $ une structure de Poisson formelle. On applique notre construction \`a ${\mathfrak g} := T_{poly}(M)[[\hbar]]$,
$ {\mathfrak g}' := {\mathcal D}_{poly}(M)[[ \hbar]]$, $ \alpha =  \hbar \gamma$ et ${\mathcal F} $ la formalit\'e de Kontsevich. 
Consid\'erons $f,g \in C^\infty(M)$. D'un c\^ot\'e, on a
\begin{eqnarray*}
 Q_\alpha ( f \otimes g)  &=  & Q_{\alpha,1} (f) \otimes g - f \otimes Q_{\alpha,1} (g) + Q_{\alpha,2}( f \otimes g) \\
 & =&  \hbar \left(  H_f \otimes g - f \otimes H_g - \{f,g\} \right).
\end{eqnarray*}
Appliquons ${\mathcal B}$ et projetons sur $D_{poly}(M)[[ \hbar]]$. Sachant que ${\mathcal B}_2( H_f \otimes g )=0 $ par d\'efinition de ${\mathcal B}_2$, on obtient la quantit\'e~:
\begin{equation}
 \label{eq:manchon1}  - \hbar  \left( {\mathcal B}_2 (f \otimes H_g)  + {\mathcal B}_1(\{f,g\} ) \right).
 \end{equation}
D'un autre c\^ot\'e, appliquons $ Q_{\beta}' \circ {\mathcal B}$ \`a  $f \otimes g $, et projetons sur $D_{poly}(M)[[  \hbar ]]$.  
Sachant que ${\mathcal B}_2(f \otimes g)=0$ pour des raisons de degr\'e, le seul terme qui reste est $ Q_{\beta,2}' ({\mathcal B}_1(f) \otimes {\mathcal B}_1(g) ) $
En comparant avec (\ref{eq:manchon1}), on obtient la formule~:
$$  - \hbar  \left( {\mathcal B}_2 (f \otimes H_g)  + {\mathcal B}_1(\{f,g\} ) \right) =  Q_{\beta,2}' ({\mathcal B}_1(f) \otimes {\mathcal B}_1(g) ) $$
Or, par d\'efinition du crochet de Gerstenhaber et du crochet d\'eriv\'e~:
$$ Q_{\beta,2}' (  {\mathcal B}_1(f),{\mathcal B}_1(g))  = [ {\mathcal B}_1(f), {\mathcal B}_1(g)]_\beta' = - [ [ \star , {\mathcal B}_1(f)]_G, {\mathcal B}_1(g)  ]_G  = - ( {\mathcal B}_1(f) \star {\mathcal B}_1(g) - {\mathcal B}_1(g) \star {\mathcal B}_1(f) ).  $$
Utilisons les notations de \cite{Manchon} et notons $ \Phi  := T^1_\alpha {\mathcal F}  $
 la d\'eriv\'ee premi\`ere de $ {\mathcal F}$ et 
$ \Psi = T^2_\alpha {\mathcal F} $  la d\'eriv\'ee seconde  en $\alpha =  \hbar \gamma$, on obtient~:
 $$  \Psi ( H_f, H_g )  = \frac{1}{\hbar} \left( \Phi(\{f,g\}) - \frac{ \Phi(f) \star \Phi(g) - \Phi(g) \star \Phi(f) }{\hbar} \right). $$
 



Nous finissons par un certain nombre de questions. Notre construction donne l'existence d'un Leib$_\infty$-morphisme entre deux alg\`ebres de Leibniz d\'eriv\'ees. Mais il serait int\'eressant de se demander ce qu'il en est des alg\`ebres de Leibniz 
g\'en\'erales. Par exemple, sous quelles conditions le th\'eor\`eme de transfert reste t-il vrai pour le contexte Leibniz~?

Par ailleurs, une question qui nous para\^it int\'eressante est de savoir ce qu'est un Maurer-Cartan pour une alg\`ebre de Leibniz diff\'erentielle gradu\'ee $ {\mathfrak g}$.
Une r\'eponse possible est de dire que c'est un \'el\'ement de type groupe $g$ et de degr\'e $0$ dans la cog\`ebre $(T({\mathfrak g}[-1]), \Delta) $ telle que $D g =0$.
Il n'est plus en g\'en\'eral vrai que les \'el\'ements de type groupe sont de la forme $e^a$ avec $a \in {\mathfrak g}_1 = {\mathfrak g}[-1]_0$.
Mais ces \'el\'ements doivent exister en g\'en\'eral. Il suit du th\'eor\`eme \ref{theo:principal} quand on a un morphisme de DGLA, un \'el\'ement de Maurer-Cartan et son image,
un morphisme d'alg\`ebres de Leibniz diff\'erentielles gradu\'ees est induit. De mani\`ere g\'en\'erale, une approche int\'eressante consiste \`a d\'eterminer les sous-ensembles stables par les
 morphisme Leibniz-infinis, les \'el\'ements de Maurer-Cartan n'\'etant qu'un exemple.

 Enfin, le cas $T_{poly}(M), D_{poly}(M)$ nous semble m\'eriter une \'etude plus approfondie. Nous avons pu d\'eriver une formule de Dominique Manchon, mais il nous semble clair que d'autres applications 
 existent.


\end{document}